\documentclass{jgg} 
\usepackage{microtype}    
\usepackage{amssymb,latexsym,amsmath,amsfonts,amsthm}
\usepackage{graphicx,psfrag,overpic,contour,color} 
\usepackage{url,cancel} 
\newtheorem{thm}{Theorem}[section]
\newtheorem{cor}{Corollary}[section] 
\newtheorem{lem}{Lemma}[section]
\theoremstyle{definition}
\newtheorem{deft}{Definition}[section]
\theoremstyle{remark}
\newtheorem{rem}{Remark}
\numberwithin{equation}{section}
\def\Square{\vbox{\hrule\hbox{\vrule height 1.8mm\hskip 1.8mm%
                         \vrule height 1.8mm}\hrule}}             
\newenvironment{Proof}{\noindent{\em Proof.}}
        {~\hfill\Square\par\medskip}

\title{The Simplest Flexible Cross-Polytopes}
\author{Hellmuth Stachel}
\affiliation{Vienna University of Technology,  
\addaddress{Wiedner Hauptstr.\ 8-10/104, 1040 Wien, Austria}{stachel@dmg.tuwien.ac.at}}
\shortauthor{H.\ Stachel}
\shorttitle{The Simplest Flexible Cross-Polytopes}

\begin{document}
\maketitle

\newcommand\bild[2]{\includegraphics[width=#2]{#1}}
\newcommand\bildh[2]{\includegraphics[height=#2]{#1}}
\newcommand\psone[2]{\centerline{\bild{#1}{#2}}}
\newcommand\psoneh[2]{\centerline{\bildh{#1}{#2}}}
\newcommand\Lemref[1]{Lemma~\ref{#1}}
\newcommand\Thmref[1]{Theorem~\ref{#1}}
\newcommand\Figref[1]{Figure~\ref{#1}}
\newcommand\Corref[1]{Corollary~\ref{#1}}
\newcommand\Defref[1]{Definition~\ref{#1}}
\def\Quadrat#1{\vbox{\hrule\hbox{\vrule height #1\hskip #1\vrule height #1}\hrule}}
\def\Beweisende{\Quadrat{2.2mm}}            
\def\BewEnde{\hfill{\Beweisende}}
\def\mBewEnde{\eqno{\Beweisende}}
\def\mklBewEnde{\eqno{\klBeweisende}}
\def\ol#1{\overline{#1}}
\def\wt#1{\widetilde{#1}}
\def\RR{{\mathbb R}}
\def\ZZ{{\mathbb Z}}
\def\NN{{\mathbb N}}
\newcommand\Vkt[1]{{\mathbf #1}}
\def\wkl{<\mskip-10mu)\mskip4mu}
\def\Frac#1#2{{\displaystyle\frac{#1}{#2}}}
\def\smFrac#1#2{\mbox{\small $\displaystyle\frac{#1}{#2}$ }}
\def\ssmFrac#1#2{\mbox{\footnotesize $\displaystyle \frac{#1}{#2}$}}
\def\sssmFrac#1#2{\mbox{\tiny $\displaystyle \frac{#1}{#2}$}}
\def\cMatrix#1{\left( \begin{array}{c} #1 \end{array}\right)}
\def\ccMatrix#1{\left( \begin{array}{cc} #1 \end{array}\right)}
\def\cccMatrix#1{\left( \begin{array}{ccc} #1 \end{array}\right)}
\def\Sum{\displaystyle\sum}
\newcommand\Ccal[1]{\mathcal C^{#1}}
\newcommand\Qcal{\mathcal Q}
\newcommand\Scal{\mathcal S}
\def\strqu{^{\prime\mskip 2mu 2}}
\def\zwstrqu{^{\prime\prime\mskip 2mu 2}}
\def\qu{^{\mskip 4mu 2}}
\def\strT{^{\prime\mskip 2mu\top}}
\def\klzwi{\mskip 1mu}
\def\zw{\mskip 2mu}
\def\Zw{\mskip 4mu}
\def\zwi{\mskip 9mu}
\def\weg{\mskip -1mu}
\def\wweg{\mskip -3mu}
\newcommand\const{\mathrm{const.}}
\def\eps{\varepsilon}
\def\phi{\varphi}

\definecolor{blaucmyk}{cmyk}{1.00,0.40,0.00,0.20} 
\def\blau{\color{blaucmyk}}
\def\blue{\color{blaucmyk}}
\definecolor{hblau}{cmyk}{0.15,0.05,0.00,0.00}
\def\hblau{\color{hblau}}
\def\magenta{\color{magenta}}
\definecolor{rotcmyk}{cmyk}{0.00,1.00,1.00,0.00}
\def\rot{\color{rotcmyk}}
\def\red{\color{rotcmyk}}
\definecolor{hhgrau}{cmyk}{0.00,0.00,0.00,0.07}
\def\hhgrau{\color{hhgrau}}
\definecolor{hgrau}{cmyk}{0.00,0.00,0.00,0.15}
\def\hgrau{\color{hgrau}}
\definecolor{grau}{cmyk}{0.00,0.00,0.00,0.50}
\def\grau{\color{grau}}
\definecolor{dgrau}{cmyk}{0.00,0.00,0.00,0.85}
\def\dgrau{\color{dgrau}}
\definecolor{gruen}{cmyk}{1.00,0.00,0.90,0.30} 
\def\gruen{\color{gruen}}
\definecolor{gelb}{cmyk}{0.00,0.00,1.00,0.00} 
\def\gelb{\color{gelb}}
\definecolor{gelbb}{cmyk}{0.00,0.00,0.50,0.00} 
\def\gelbb{\color{gelbb}}
\definecolor{gelbbb}{cmyk}{0.00,0.00,0.20,0.00} 
\def\gelbbb{\color{gelbbb}}
\def\white{\color{white}}
\definecolor{weiss}{gray}{1}
\newcommand{\weiss}{\color{weiss}}
\definecolor{lila}{rgb}{0.8,0.0,0.8}
\newcommand{\violett}{\color{lila}}
\def\sz{\small}
\def\Raum#1{{\mathbb #1}}
\newcommand\euklRaum[1]{\Raum E^{#1}}

\def\Bq{B'}
\def\Bqq{B''}
\def\Sq{S'}
\def\Sqn{S^{\prime\mskip1mu s}}
\def\Sqq{S''}
\def\Bqn{B^{\prime\mskip1mu s}}
\def\Bqqn{B^{\prime\prime \mskip1mu s}}
\def\tupel#1#2#3{#1_{#2},\ldots,#1_{#3}}
\newcommand\plus{\!+\!}
\newcommand\minus{\!-\!}
\newcommand\mal{\!\times\!}
\def\dist#1#2{\ol{#1#2}}
\def\mathtext#1{\mskip 12mu \mbox{#1}\mskip 12mu}
\def\inv{^{-1}}
\def\span#1#2{[#1,#2]}
\def\adj{^\ast}
\def\poly#1{{\mathcal #1}}
\def\qpoly#1{\ol{\poly #1}}
\newcommand\Wing{\Delta}
\newcommand\CoWing{\Gamma\weg}
\newcommand\Gang{\Wing_0}
\newcommand\Gangn{\Gamma_{\wweg m}}
\newcommand\Gangq{\Gang'}
\newcommand\Gangqq{\Gang''}
\newcommand\Gangnq{\Gangn'}
\newcommand\Gangnqq{\Gangn''}
\newcommand\Ganghe{H_{\weg m}} 
\newcommand\Rast{\CoWing_0}
\newcommand\Rasthe{H_{\wweg f}} 
\newcommand\Rastn{\Gamma_{\mskip -4mu f}}
\newcommand\Achse{h} 
\newcommand\Kante{v} 
\newcommand\BKante{w} 
\newcommand\Stroph{\mathcal S}
\newcommand\comment[1]{{\violett #1}}

\hfill{Draft \today} 

\MSC[52C25]{52B11, 53A17, 70B15}


\begin{abstract}
According to R.\ Bricard there exist three types of flexible octahedra. 
The octahedra of Type~3 are unsymmetric and admit two flat poses.
With regard to higher-dimensional analogues of octahedra called cross-polytopes, A.A.\ Ga\u\i fullin presented in 2014 a complete classification of flexible types in $n$-dimensional Euclidean, hyperbolic and spherical spaces for $n>3$.
The goal of this presentation is a synthetic approach to a particular family in the Euclidean $n$-space, the flexible cross-polytopes that admit two poses within hyperplanes.
We provide a construction of their flat poses and prove several properties of these higher-dimensional analogues to Bricard's type-3 octahedra.
According to Ga\u\i fullin, they are the simplest from the algebraic point of view.
\end{abstract}

\section{Introduction}\label{sec:Intro}
In 1897, Raoul Bricard \cite{Bricard1} proved that in 3-dimensional Euclidean space $\Raum E^3$ there exist three types of flexible octahedra, i.e., polyhedra of the combinatorial type of a regular octahedron with rigid faces and edges functioning as hinges, while self-intersections are ignored.
In 1912, G.T.\ Bennett \cite{Bennett} studied the flexions of these types with regard to their spherical images. 
A different approach is due to R.\ Connelly in 1978 \cite{Connelly}. 
The author classified in \cite{Sta_35} the flexible octahedra due to their relation to confocal quadrics and based on a configuration theorem on bipartite graphs.
Flexible octahedra can also be seen as Kokotsakis meshes with a triangular basis (note, e.g., \cite{Sta_128}).
The kinematics of these meshes was recently studied in \cite{Hu}.
In \cite{Guest}, flexible open chains of type-3 octahedra without link interference were presented, where consecutive octahedra share an edge and the adjacent spanned planes.

\begin{figure}[htb]  
  \psone{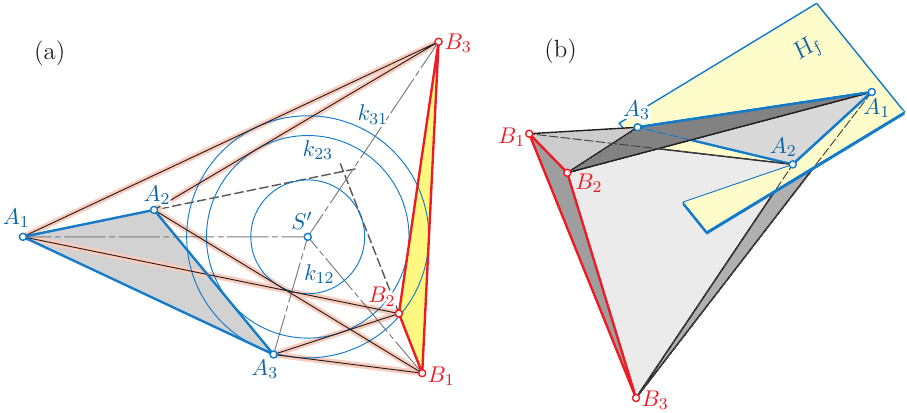}{153mm} 
  \caption{(a) Flat pose of a flexible type-3 octahedron: 
Two concentric circles $k_{12}$, $k_{31}$ and two points $A_1$, $B_1$ define the remaining vertices $A_2,\dots,B_3$ and a third contacting circle $k_{23}$. \
 (b) A spatial pose of this octahedron with the fixed triangle $A_1A_2A_3$ in the plane $\Rasthe$ (yellow).
At the depicted physical model, the two opposite faces $B_1A_2A_3$ and $A_1B_2B_3$ need to be omitted in order to avoid interpenetrations.}
  \label{fig:Fig_1}
\end{figure}

Bricard's octahedra of Type~3 are unsymmetric and doubly collapsible, i.e., they admit two flat poses.
From the geometric point of view, this type is the most complicated.
However, from the algebraic point it is the simplest as A.A.\ Ga\u\i fullin mentions in \cite[p.~90]{Gaifullin}. 
\Figref{fig:Fig_1}a shows how a flat pose of a type-3 octahedron can be constructed according to \cite[Fig.~8]{Bennett} and \cite[Fig.~297, p.~330]{Bricard2}.

\medskip
The $n$-dimensional analogues of octahedra are called {\em cross-polytopes} $\Ccal n$.
Following the notation in \cite{Gaifullin}, they have $2n$ vertices coupled into pairs of {\em opposite} vertices $(A_i, B_i)$ for $i=1,\dots, n$.
The $2^n$ hyperfaces (or {\em facets}) of $\Ccal n$ are the simplices $X_1\dots X_n$ where $X_i$ stands either for $A_i$ or $B_i$.
The $4 {n\choose 2} = 2n(n-1)$ edges of $\Ccal n$ are $X_iX_j$ for $i\ne j$. 
More general, the $k$-faces, $1 < k\le n-1$ of $\Ccal n$ are the $2^{k+1}{n\choose {k+1}}$ simplices with $k+1$ vertices $X_{i_0}\dots X_{i_k}$ out of $\{A_1,\dots,A_n, B_1,\dots,B_n\}$ such that the indices $i_0\dots i_k$ are mutually different.
Hence, no $k$-face contains a pair $(A_i,B_i)$ of opposite vertices.
We presuppose that no $k$-face is of dimension smaller than $k$, i.e., flat.

In 2014, the young Russian mathematician Alexander A.\ Ga\u\i fullin \cite{Gaifullin} surprised the scientific community with the complete solution of a long-lasting open problem, namely the question for flexible cross-polytopes in the $n$-dimensional Euclidean, hyperbolic and spherical spaces for $n > 3$.
So far, only particular flexible examples in $\Raum E^4$ where known (see \cite{Sta_88}).
Based on algebraic methods, Ga\u\i fullin succeeded to classify the flexible types in all spaces, and he even presented parametrizations of the flexions in terms of Jacobian elliptic functions.
In several papers Ga\u\i fullin addressed the problem whether the volumes of flexible cross-polytopes in non-Euclidean spaces remain constant during the flexions (see, e.g., \cite{Gaifullin2, Gaifullin3}). 

According to \cite{Gaifullin}, $n$-dimensional polytopes are called {\em flexible}, if they have hinges at their edges and the polytope can deform while each $(n\minus 1)$-dimensional facet remains congruent to itself. 
In other words, only the dihedral angles between two neighboring facets can vary.
We are interested in nontrivial deformations, which means that the deformations are not induced by a motion of a rigid cross-polytope in the ambient space.

\begin{figure}[t] 
  \psone{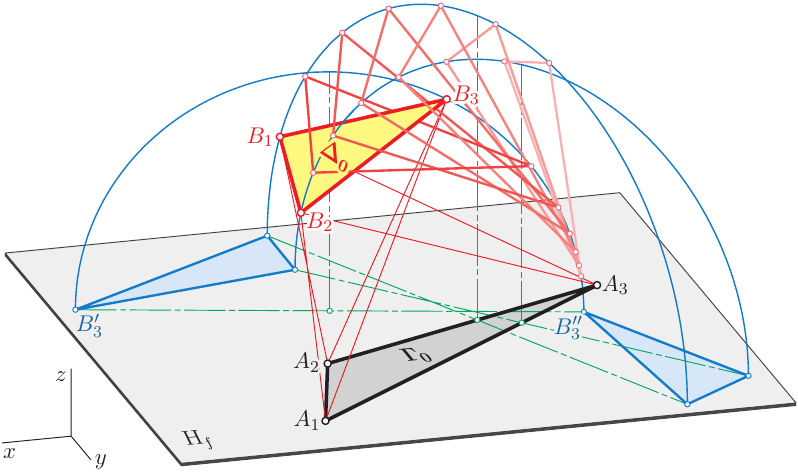}{135mm} 
  \caption{A type-3 octahedron in motion, while the triangle $\Rast = A_1A_2A_3$ remains fixed in the plane $\Rasthe\!:\, z=0$.}
  \label{fig:in_motion}
\end{figure}

When studying the flexibility of cross-polytopes $\Ccal n$, we follow the convention in \cite{Gaifullin} and assume that the simplex $\Rastn:= A_1\dots A_n$ in the hyperplane $\Rasthe$ is fixed while the opposite simplex $\Gangn:= B_1{\dots}B_n$ in the hyperplane $\Ganghe$ is moving.
This means that the vertex $B_i$ traces a circle during its rotation about the $(n\minus 2)$-dimensional {\em axis} $\Achse_i\subset\Rasthe$ spanned by $A_1, \dots, A_{i-1}, A_{i+1}, \dots, A_n$.
Point $B_i$ rotates together with the simplex called {\em wing} $\Wing_i:= A_1 \dots A_{i-1} B_i A_{i+1} \dots A_n$ relative to the fixed simplex $\Rastn$ (note the flexible octahedron in \Figref{fig:in_motion} as an example for $n=3$).

At flexible cross-polytopes, there exist simultaneous rotations of all $n$ wings $\Wing_1,\dots,\Wing_n$ such that all ${n\choose 2}$ edges of the moving simplex $\Gangn$ preserve their lengths $\ol{B_iB_j}$.\footnote{
Ga\u\i fullin created the pretty name {\em butterfly} for the union of the fixed simplex with all wings.}
Since the moving simplex together with the wings contain all edges of $\Ccal n$, every other facet remains rigid, too.
We call such a continuous self-motion of $\Ccal n$ a {\em flexion} and ignore mutual intersections between facets.

\medskip
Below, we study in the Euclidean $n$-space $\Raum E^n$ those particular flexible cross-polytopes that pass during their flexion through two flat poses, i.e., poses where the moving simplex $\Gangn$ lies in the hyperplane $\Rasthe$ spanned by the fixed simplex $\Rastn$.

\begin{deft}\label{def:twice-flat}
A cross-polytope $\Ccal n$ is called {\em twice-flat} if it admits two flat poses such that for each vertex $B_i$ of the moving simplex the two corresponding positions $\Bq_i$ and $\Bqq_i$ are different.
\end{deft}

According to \cite[eq.\,(3.3)]{Gaifullin}, the constancy of the edge length $\ol{B_iB_j}$, $i\ne j$, during the flexion is equivalent to a particular biquadratic relation between the tangents $t_i$ and $t_j$ of halved dihedral angles $\phi_i,\phi_j$ between $\Rastn$ and the respective wings $\Wing_i, \Wing_j$ (see also \cite[eq.\,(4)]{Sta_128}).
Since at flexible twice-flat cross-polytopes the values $t_i = t_j = 0$ and $t_i = t_j = \infty$ must satisfy these conditions for each pair $(i,j)$ with $i,j\in\{0,\dots,n\}$ and  $i\ne j$, all biquadratic relations split.
In this case for each pair of wings $\Wing_i, \Wing_j$ two simultaneous one-parameter rotations about the respective axes $\Achse_i, \Achse_j$ through angles with proportional $t_i$ and $t_j$ are available for preserving the distance $\ol{B_iB_j}$.
While Ga\u\i fullin showed the existence of flexible versions in an algebraic way, the goal of our study is a synthetic approach, a geometric characterization of the flat poses and a geometric analysis of the general poses of flexible twice-flat cross-polytopes in the Euclidean $n$-space. 

\medskip
The paper is structured as follows:
In the coming section we revisit Type~3 of Bricard's octahedra.
Our novel and mainly synthetic approach yields for the flat poses of twice-flat octahedra two geometric conditions that are necessary and sufficient for flexibility.
In addition, we disclose several properties of {\em spatial}, i.e., non-flat poses, namely symmetries between planes spanned by the faces. 

Section~3 is devoted to $n$-dimensional flexible twice-flat cross-polytopes $\Ccal n$ in $\Raum E^n$.
Again, we come up with conditions that characterize flat poses of flexible versions.
A construction of flexible twice-flat cross-polytopes is given that generalizes Bricard's construction in \cite[p.\,144]{Bricard1}.    
Finally, we check which properties of spatial poses at type-3 octahedra can be generalized to higher dimensions.

\section{Revisiting Bricard's type-3 octahedra}
Let $\Ccal 3$ be a cross-polytope in $\Raum E^3$, i.e., an octahedron with the triangles  $\Rast:= A_1A_2A_3$ and $\Gang:= B_1B_2B_3$ that serve respectively as fixed and moving simplex $\Rastn$ and $\Gangn$.\footnote{
Because of particular symmetries (note \Figref{fig:spher_linkage2}) it makes sense to use in this chapter the new symbols $\Rast$ and $\Gang$ instead of $\Rastn$ and $\Gangn$ for the fixed and the moving simplex.}
When $\Ccal 3$ happens to admit two different flat poses, then in these poses the vertices $B_1$, $B_2$ and $B_3$ of $\Gang$ will be marked by one or two primes for better distinction from spatial poses.

\subsection{Local symmetries of twice-flat octahedra}

\begin{thm}\label{thm:symmetry3} 
A flat octahedron with the fixed triangle $A_1A_2A_3$ and the coplanar opposite triangle $\Bq_1\Bq_2\Bq_3$ admits a second flat pose $\Bqq_1\Bqq_2\Bqq_3$ with $\Bq_i\ne \Bqq_i$ for each $i\in\{1,2,3\}$ if and only if at each vertex $A_i$ the connecting lines with the remaining two pairs $(A_j,\Bq_j)$, $j\ne i$, of opposite vertices have common axes of symmetry.
\end{thm}

\begin{figure}[htb]  
  \psone{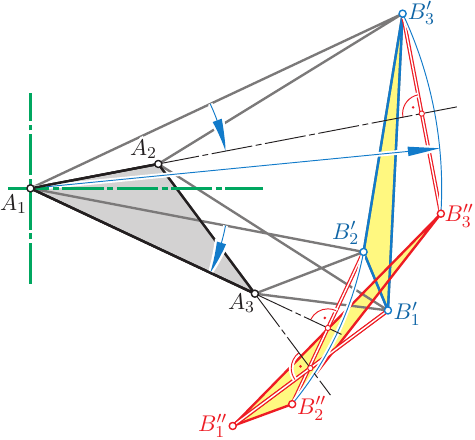}{80mm}
  \caption{At a twice-flat octahedron the moving triangle $\Gang$ varies between $\Bq_1\Bq_2\Bq_3$ and $\Bqq_1\Bqq_2\Bqq_3$ relative to the fixed triangle $\Rast = A_1A_2A_3$.
The double lines are projections of the circular paths traced by the moving vertices in the case of a flexible octahedron.}
  \label{fig:Fig_2}
\end{figure}

\begin{Proof}
Suppose that $\Ccal 3$ admits two flat poses with $\Bq_i\ne \Bqq_i$ for all $i$.
Then the halfturn in $\Raum E^3$ about the axis $[A_1, A_2]$ sends $\Bq_3$ to $\Bqq_3$, and that about $[A_1, A_3]$ sends $\Bq_2$ to $\Bqq_2$ (see \Figref{fig:Fig_2}).\footnote{
~Throughout the paper we use brackets $[X,Y,\dots,Z]$ as the symbol for the affine hull (or span) of the listed points.
For example, $[X,Y]$ stands for the {\em line} connecting the points $X$ and $Y$, and
$[\Wing_i]$ stands for the hyperplane spanned by the wing $\Wing_i$.
On the other hand, the symbol $XY\dots Z$ denotes the simplex with vertices $X,Y,\dots,Z$, and in particular $XY$ denotes the {\em segment} terminated by $X$ and $Y$.} 
Thus, within the plane $\Rasthe = [A_1,A_2,A_3]$ the bisectors of the point pairs $(\Bq_2, \Bqq_2)$ and $(\Bq_3, \Bqq_3)$ intersect at $A_1$, while there are equal distances $\ol{\Bq_2\Bq_3} = \ol{\Bqq_2\Bqq_3}$.
This implies that $A_1$ is the center of a rotation within $\Rasthe$ that brings $\Bq_2\Bq_3$ to $\Bqq_2\Bqq_3$.
One half of the angle of rotation is marked in \Figref{fig:Fig_2} which reveals that the pairs of lines $([A_1, A_2], [A_1, \Bq_2])$ and $([A_1, \Bq_3], [A_1, A_3])$ share the angle bisectors. 

\noindent
The disclosed symmetry at $A_1$ occurs similarly at $A_2$ and $A_3$, which confirms the claim.
\end{Proof}

For the existence of a second flat pose it is not relevant which face of the octahedron is fixed.
This implies

\begin{cor}\label{cor:symmetry3} 
Referring to \Thmref{thm:symmetry3}, the stated symmetries at $A_1$, $A_2$ and $A_3$ are equivalent to common axes of symmetries at each point $\Bq_i$, $i\in\{1,2,3\}$,  between the connections with the remaining two pairs $(A_j,\Bq_j)$ with $j\ne i$. 
\end{cor}

\begin{deft}\label{def:loc_symm 2d} 
A flat pose $A_1\dots \Bq_3$ of an octahedron is called {\em locally symmetric}, when at each vertex $A_i$ and $\Bq_i$ the respective connections with the remaining two pairs $(A_j,\Bq_j)$ for $j\ne i$, have common axes of symmetry.
\end{deft}

It will turn out that local symmetry of a flat octahedron is not sufficient for the flexibility of a twice-flat octahedron. 
Below we summarize a few properties of a locally symmetric pose.

\begin{figure}[htb]  
  \psone{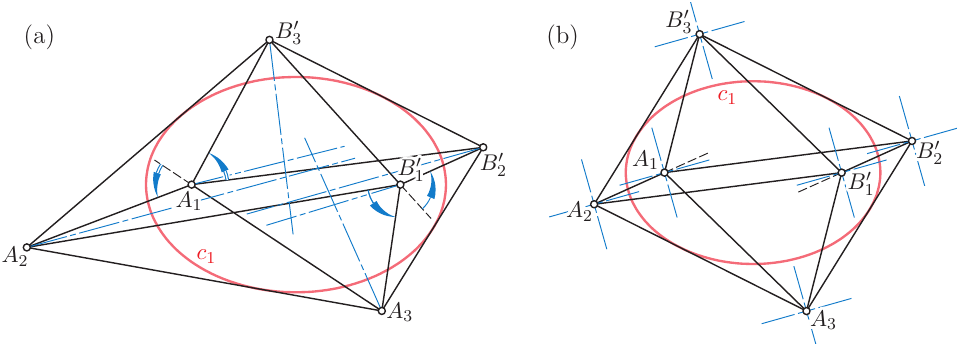}{162mm}
  \caption{(a) A flat pose $A_1\dots \Bq_3$ is locally symmetric if and only if $A_1$ and $\Bq_1$ are focal points of a conic $c_1$ that contacts the sides of the quadrangle $A_2A_3\Bq_2\Bq_3$ (\Lemref{lem:inscribed_conic}).\\
(b) A locally symmetric flat pose with parallel axes of symmetry between $[A_i,A_j]$ and $[A_i,B_j]$ as well as between $[B_i,A_j]$ and $[B_i,B_j]$ for all $i\ne j$. 
In this case all quadrangles $A_j\dots \Bq_k$ are parallelograms.}
 \label{fig:twoparam_Ks_schar}  \label{fig:parallel_axes}
\end{figure}

\begin{lem}\label{lem:inscribed_conic} 
{\em (i)} A flat pose $A_1\dots \Bq_3$ is locally symmetric if and only if $A_1$ and $\Bq_1$ are focal points of a conic $c_1$ that contacts the four sides of the quadrangle $A_2A_3\Bq_2\Bq_3$ (\Figref{fig:twoparam_Ks_schar}a). 
\\[0.5mm]
{\em (ii)} At a locally symmetric flat pose the eight pedal points of $A_i$ and $\Bq_i$ on the sides of the quadrangle $A_jA_k\Bq_j\Bq_k$ are concircular for each permutation $(i,j,k)$ of $(1,2,3)$.
Moreover, the midpoints of all pairs $(A_i,\Bq_i)$ are collinear.
\\[0.5mm]
{\em (iii)} At twice-flat octahedra, the two flat poses $\Bq_1\Bq_2\Bq_3$ and $\Bqq_1\Bqq_2\Bqq_3$ of the moving triangle are indirectly congruent, and the midpoints of $(\Bq_i,\Bqq_i)$, $i=1,2,3$, are collinear.
\end{lem}

\begin{Proof}
(i) If a locally symmetric flat pose $A_1\dots\Bq_3$ is given, then due to Desargues's involution theorem (see, e.g., \cite[p.~336]{Conics}) the connections of $A_1$ with the pairs $(A_2,\Bq_2)$ and $(A_3,\Bq_3)$ define a symmetric involution.
Since its fixed lines are orthogonal, the isotropic lines through $A_1$ correspond each other in this involution.
\\
Consequently, there exists a conic $c_1$ with focal point $A_1$ that contacts the four sides of the quadrangle $A_2A_3\Bq_2\Bq_3$ (note, e.g., \cite[p.~288]{Conics}).
Since the tangents from any point $X$ to a conic and the connections with the focal points share the angle bisectors (see, e.g., \cite[p.~43]{Conics}), the symmetries at the vertices $A_2,A_3,\Bq_2,\Bq_3$ of the contacting quadrangle imply that the point $B_1$ is the second focal point.
\\[0.5mm]
Conversely, let a conic $c_1$ with focal points $A_1,\Bq_1$ and two different points $A_2,\Bq_2$ in the exterior of $c_1$ be given.
We specify $(A_3,\Bq_3)$ as another pair of opposite points of the quadrilateral formed by the tangents drawn from $A_2$ and $\Bq_2$ to $c_1$.
The common angle bisectors between the tangents from any point $X$ to $c_1$ and the connections with the focal points imply the required symmetries at $A_2, \dots, \Bq_3$. 
\\
The symmetries at $A_1$ and $\Bq_1$ follow from Desargues's involution theorem.
The tangents from $A_1$ or $B_1$ to all curves of the range defined by the sides of the quadrangle $A_2\dots\Bq_3$ belong to an involution, in particular the connections with $(A_2,\Bq_2)$ as well as that with $(A_3,\Bq_3)$. 

\smallskip\noindent
(ii) The pedal points of the sides of the quadrangle $A_j\dots\Bq_k$ with respect to (w.r.t.\ in brief) the focal points $A_i$ and $\Bq_i$ are located on a circle (see \cite[Fig.~2.16]{Conics}).
All conics tangent to the sides of $A_jA_k\Bq_j\Bq_j$ belong to a range (dual pencil).
Therefore the conics' centers are located on a line, that passes also through the midpoints of $(A_j,\Bq_j)$ and $(A_k,\Bq_k)$ as the carriers of pairs of line pencils contained in the range. 

\smallskip\noindent
(iii) The congruence transformation that takes $\Bq_1\Bq_2\Bq_3$ to $\Bqq_1\Bqq_2\Bqq_3$ must be orientation reversing since otherwise it would simultaneously be equal to the three rotations with respective centers $A_i$, $i=1,2,3$, that send the single side $\Bq_j\Bq_k$ to its respective image $\Bqq_j\Bqq_k$.
Thus, it is a glide reflection along an axis that passes through the midpoints of corresponding points.
\end{Proof}

Statement (iii) offers another way for obtaining a locally symmetric flat pose:
Choose two indirectly congruent triangles $\Bq_1\Bq_2\Bq_3$ and $\Bqq_1\Bqq_2\Bqq_3$ in a common plane.
If the bisectors of $(\Bq_i,\Bqq_i)$ for $i=1,2,3$ form a triangle, then this serves as the fixed triangle $A_1A_2A_3$.

\medskip
The flat poses of twice-flat octahedra in \Figref{fig:Fig_2} show that at the four-sided pyramid connecting $A_1$ with the sides of the quadrangle $A_2A_3B_2B_3$ opposite interior angles at $A_1$ are congruent (also called {\em Voss condition} \cite{Kilian}). 
Similar congruences appear at all other vertices.
However, when the flat pose of the pyramid is not twofold covered (see $A_2$ in \Figref{fig:Fig_1}a), then opposite angles sum up to $\pi$ ({\em anti-Voss condition} by \cite{Kilian}).

\begin{deft} 
A four-sided pyramid is called {\em isogonal} if either opposite interior angles at the apex are congruent or they are supplementary, i.e., they sum up to $\pi$.
The extended lines form an {\em isogonal double-pyramid} consisting of four concurrent lines $g_1,\dots, g_4$ in cyclic order where opposite angles are equal, i.e., $\wkl g_1g_2 = \,\wkl g_3g_4$ and $\wkl g_2g_3 = \,\wkl g_4g_1$ under $0 \le \,\wkl g_ig_j \le \frac\pi 2\zw$.
\end{deft}

\begin{figure}[htb] 
  \psone{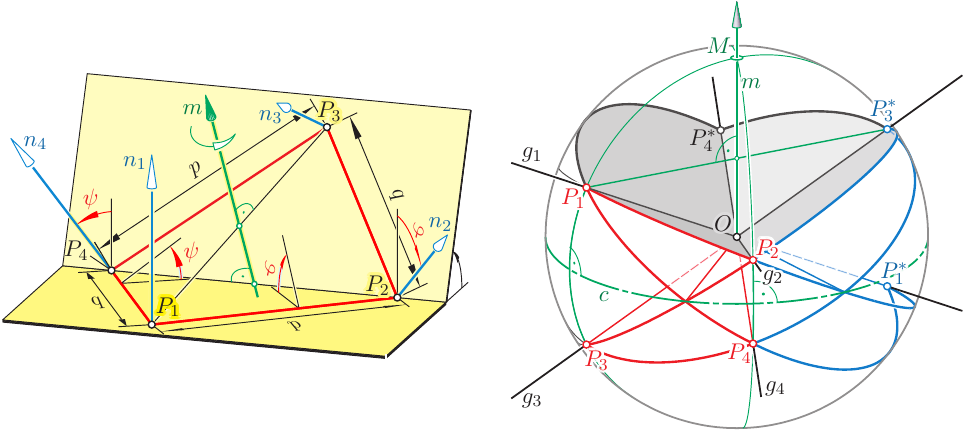}{163mm}
  \caption{Left: Each skew isogram $P_1P_2P_3P_4$ has an axis of symmetry $m$, which
meets both diagonals orthogonally at their midpoints. 
\\
Right: At the depicted isogonal double-pyramid with apex $O$ and generators $g_1,\dots,g_4$ the three isogonal pyramids through $P_1P_2P_3P_4$, $P_1^\ast P_2 P_3^\ast P_4$ and $P_1P_2P_3^\ast P_4^\ast$ are emphasized. 
Two pyramids are crossed and symmetric w.r.t.\ to the plane $\mu$ of the circle $c$.
The non-crossed pyramid is symmetric w.r.t.\ the circle's axis $m$.
At each pyramid opposite dihedral angles are congruent.}
  \label{fig:double-cone}
\end{figure}

An isogonal four-sided pyramid intersects a sphere centered at the apex along an {\em isogram}, i.e., a quadrangle with opposite sides of equal length (\Figref{fig:double-cone},\,left).
In our case, we obtain a spherical isogram where the straight sides are replaced by arcs on great circles.
A spherical isogram is either free of self-intersections or two opposite sides are intersecting quite similar to a planar antiparallelogram.
In the latter case we speak of a {\em crossed} spherical isogram (\Figref{fig:double-cone},\,right).
By replacing one edge of the pyramid by its complement, i.e., by the halfline pointing in the opposite direction, supplementary opposite angles become congruent angles, and vice versa. 
When two adjacent edges of an isogonal pyramid are replaced by their complements, then a crossed isogram is replaced by an uncrossed version, and vice versa (note the isogonal double-pyramid in \Figref{fig:double-cone},\,right, containing three pyramids with congruent opposite angles).

\goodbreak
\begin{lem}\label{lem:symm_doublePyr} 
If the four lines $g_1, \dots, g_4$ of an isogonal double-pyramid are not coplanar, then the line $m$ of intersection of the two diagonal planes $[g_1,g_3]\cap [g_2,g_4]$ is an axis of symmetry.
Dually, the plane $\mu$ connecting $[g_1,g_2]\cap[g_3,g_4]$ with $[g_2,g_3]\cap[g_4,g_1]$ is orthogonal to $m$ and a plane of symmetry. 
\end{lem}

\begin{Proof}
We select on the double-pyramid a pyramid with congruent opposite interior angles.
Then the resulting isogonal pyramid intersects the unit sphere centered at the apex $O$ along a skew spherical isogram (Figures~\ref{fig:double-cone},\,right, or \ref{fig:sph_isogram}).
For each spherical isogram holds that the common perpendicular of the diagonals passes through their midpoints and is the axis of a halfturn with spherical center $M$ that exchanges opposite vertices (see, e.g., \cite[p.~555]{Quadrics}).
This axis passes also through the sphere's center $O$ because of equal distances to the vertices of the isogram.
The halfturn about the diameter $[O,M]$ exchanges opposite lines of the original double-pyramid.
\\
After composition with the reflection in the center $O$, the reflection in the diameter plane perpendicular to $[O,M]$ sends the double-pyramid onto itself, too. 
This diameter plane $\mu$ contains the intersections between corresponding planes $[g_1,g_2]\cap[g_3,g_4]$ and $[g_1,g_4]\cap[g_2,g_4]$. 
In \Figref{fig:double-cone},\,right, the great circle $c$ serves on the sphere as the axis of symmetry of a crossed spherical isogram.
\end{Proof}

\smallskip
Four-sided pyramids (without base quadrangle) are flexible (\Figref{fig:sph_isogram}).
It is wellknown (see, e.g., \cite[eq.~(9)]{Sta_128}) that during the flexion of an isogonal pyramid the tangents of the halved dihedral angles $\phi_{3}$ along $A_1A_2$ and $\phi_{2}$ along $A_1A_3$ are proportional (note \Lemref{lem:pyramid bending} below).

\subsection{Transmission between adjacent wings}

\begin{figure}[htb]  
  \psone{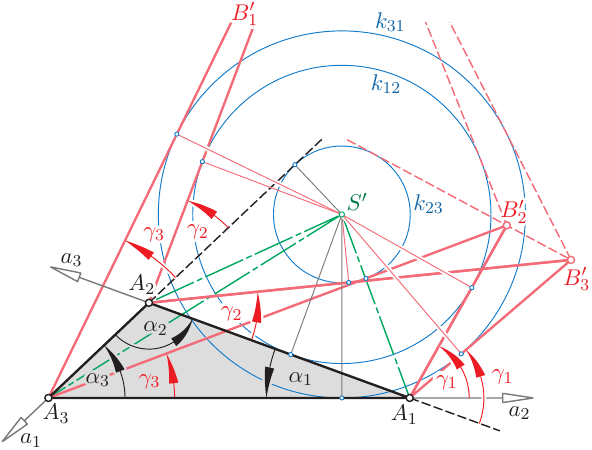}{100mm}
  \caption{The angles in the fixed face $\Rast = A_1A_2A_3$ and in the wings $\Wing_1 = A_2A_3B_1$, $\Wing_2 = A_3A_1B_2$, and $\Wing_3 = A_1A_2B_3$ of a twice flat octahedron that is displayed in a flat pose.}
  \label{fig:denotation}
\end{figure}

Let $A_1\dots\Bq_3$ be a locally symmetric flat pose of an octahedron.
Suppose that the octahedron is flexible.
For obtaining an exact formula for the transmission between the wings, we assume that the fixed triangle $\CoWing_0 = A_1A_2A_3$ is counter-clockwise oriented, and the directions of the axes of rotations $a_1, a_2, a_3$ are given, respectively, by the directed segments $A_2A_3$, $A_3A_1$ and $A_1A_2$.
The signed interior angles of $A_1A_2A_3$ are denoted by $\alpha_1, \alpha_2, \alpha_3$ (\Figref{fig:denotation}).

With $\gamma_1$ we denote the signed angle at $A_1$ between the lines $[A_1,A_3]$ and $[A_1,\Bq_2]$.
It equals the angle between $[A_1,A_2]$ and $[A_1,\Bq_3]$.
We use the rotation by $\gamma_1$ about $A_1$ to transfer the orientation of $a_2$ to the line $[A_1,\Bq_2]$ and denote with $p_1$ the signed distance of $A_1\Bq_2$, and similarly with $q_1$ that of $A_1\Bq_3$.
An analogous procedure yields the signed angles $\gamma_2, \gamma_3$ and the signed distances $p_2, q_2$ and $p_3,q_3$ at the remaining vertices $A_2$ and $A_3$.
We assume for all $i$ that $0 < \alpha_i,\gamma_i < \pi$ and $(\alpha_i+\gamma_i), (\alpha_i-\gamma_i)\ne k\pi$ for all $k\in\mathbb Z$.

\begin{lem}\label{lem:pyramid bending} 
Referring to the previous notation, the angles of rotation $\phi_2$ of the wing $\Wing_2$ about $a_2$ and $\phi_3$ of $\Wing_3$ about $a_3$ relative to $\Rast$ with $t_i:= \tan\frac{\phi_i}2$ for $i=2,3$ are related by
\[  t_3 = \frac{-\sin\gamma_1 \pm \sin\alpha_1}{\sin(\gamma_1 + \alpha_1)}\,t_2\,,
   \quad\mbox{hence} \quad 
    t_3 = \frac{\sin\frac{\alpha_1 - \gamma_1}2}{\sin\frac{\alpha_1 + \gamma_1}2}\,t_2 
    \quad \mbox{or} \quad     
    t_3 = -\frac{\cos\frac{\alpha_1 - \gamma_1}2}{\cos\frac{\alpha_1 + \gamma_1}2}\,t_2
\]
according to the upper or the lower sign in the first equation.
\end{lem}

\begin{figure}[hbt] 
  \psone{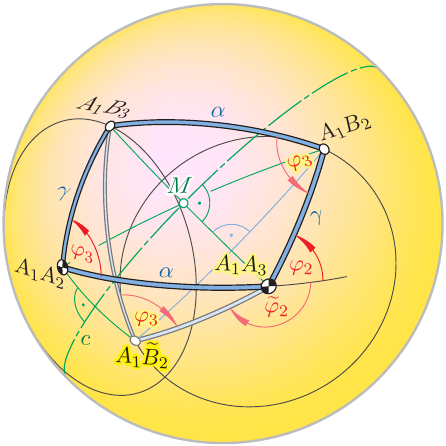}{75mm}
  \caption{The flexion of the isogonal pyramid with apex $A_1$ induces on a sphere centered at $A_1$ two rational spherical motions with bifurcations at the aligned poses.
A crossed spherical isogram is met by its spherical axis $c$ of symmetry; otherwise the isogram includes its center $M$ of symmetry.}
\label{fig:sph_isogram} 
\end{figure}

\begin{figure}[htb] 
  \psone{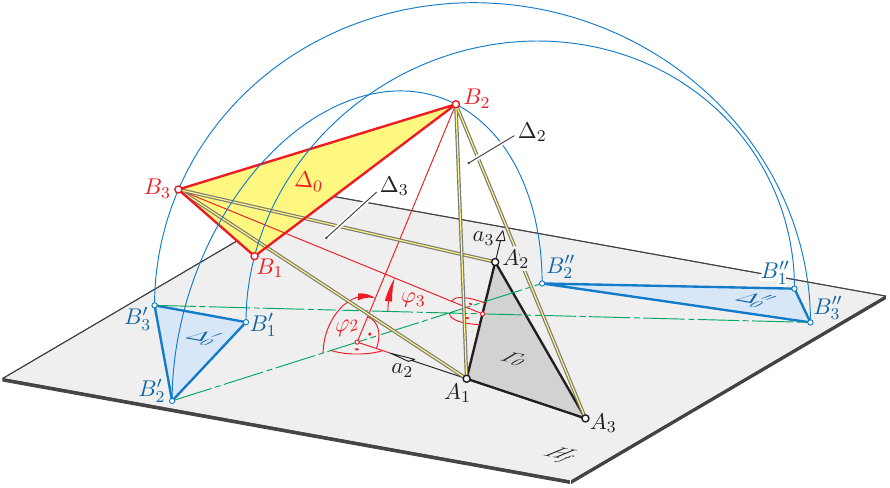}{150mm}
  \caption{Spatial pose of the moving triangle $B_1B_2B_3$ with the corresponding angles of rotation $\phi_2$, $\phi_3$ of $B_2, B_3$ and wings $\Wing_2 = A_1A_3B_2$ and $\Wing_3 = A_1A_2B_3$ along with the respective axes $a_2$ and $a_3$.}
 \label{fig:moving}
\end{figure}

\begin{Proof}
We use a right-handed coordinate frame with the flat pose in the plane $z=0$, with $A_1$ as origin and $a_2$ as positive $x$-axis.
From points in the half-space $z>0$ the triangle $A_1A_2A_3$ appears counter-clockwise oriented, as shown in \Figref{fig:denotation}.

Thus, the circular paths of $B_2$ and $B_3$ (\Figref{fig:moving}) can be parametrized as
\setlength\arraycolsep{1.5mm}
\[  B_2 = p \cMatrix{\cos\gamma_1\\ \sin\gamma_1\cos\phi_2\\ \sin\gamma_1\sin\phi_2}\weg, \
    B_3 = q \cMatrix{-\cos\alpha_1 \cos\gamma_1 - \sin\alpha \sin\gamma_1 \cos\phi_3\\
                \sin\alpha_1 \cos\gamma_1 - \cos\alpha \sin\gamma_1 \cos\phi_3\\
                \sin\gamma_1 \sin\phi_3}\weg .   
\]  
After some computations (for details see, e.g., \cite{Sta_128}) the claim that the distance $\ol{B_2B_3}$ remains equal to $\ol{\Bq_2\Bq_3}$ in the initial flat pose with $\phi_2 = \phi_3 = 0$ is equivalent to
\[  \sin(\alpha_1 - \gamma_1)\,t_2^2 - \sin(\alpha_1 + \gamma_1)\,t_3^2 - 2\sin\gamma_1\, t_2 t_3 = 0.
\]
This reveals that in the particular case of an isogonal pyramid the biquadratic relation (according to \cite[eq.~(4)]{Sta_128} or \cite[eq.~(3.3)]{Gaifullin}) splits into two linear functions
\[  t_3 = \frac{\weg -\sin\gamma \pm \sqrt{\sin^2\weg\gamma + \sin^2\weg\alpha \cos^2\weg\gamma - \cos^2\weg\alpha \sin^2\weg\gamma}}{\sin(\alpha + \gamma)}\,t_2 
    = \frac{\weg -\sin\gamma \pm \sin\alpha}{\sin(\alpha + \gamma)}\,t_2.
\]
Thus we obtain either
\begin{equation}\label{eq:transmission0}
  t_3 = \frac{\sin\frac{\alpha - \gamma}2}{\sin\frac{\alpha + \gamma}2}\,t_2
   \quad \mbox{or} \quad
    t_3 = -\frac{\cos\frac{\alpha - \gamma}2}{\cos\frac{\alpha + \gamma}2}\,t_2\,.
\end{equation}
Both flat poses allow bifurcations between the two rational movements. 
\Figref{fig:sph_isogram} shows the spherical images of two poses, one of each movement.
\end{Proof}

\begin{lem}\label{lem:ratio} 
Let $A_1\dots \Bq_3$ be a locally symmetric flat pose of an octahedron.
Then, referring to \Lemref{lem:pyramid bending}, all points $S\ne A_1$ of one axis of  symmetry between the lines $[A_1,A_2], [A_1,\Bq_2]$ and $[A_1,A_3], [A_1,\Bq_3]$ satisfy the relation	
\[  t_3 = \frac{d(S,a_2)}{d(S,a_3)}\,t_2\,,
\]
where $d(S,a_i)$ for $i=1,2,3$ denotes the signed distance of $S$ to the revolute axis $a_i$ such that for points in the halfplane of $a_i$ enclosing the interior of $A_1A_2A_3$ the distance is positive.
In other words, $\left(d(S,a_1),\,d(S,a_2),\,d(S,a_3)\right)$ are the normalized trilinear coordinates of $\Sq$ w.r.t.\ the fixed triangle $A_1A_2A_3$. 
\end{lem}

\begin{Proof}
Let $(i,j,k)$ by a cyclic permutation of $(1,2,3)$.
If one symmetry axis through $A_i$ meets the opposite side $A_jA_k$, then it includes with the axes $a_k$ and $a_j$ the angles $\ssmFrac{\alpha_i-\gamma_i}2$ and $\ssmFrac{\alpha_i+\gamma_i}2$, respectively (note vertex $A_3$ in \Figref{fig:denotation}).
Therefore holds, by virtue of \Lemref{lem:pyramid bending},
\[  \frac{-\sin\gamma_i + \sin\alpha_i}{\sin(\gamma_i + \alpha_i)} 
   = \frac{\sin\ssmFrac{\alpha_i - \gamma_i}2}{\sin\ssmFrac{\alpha_i + \gamma_i}2} 
   = \frac{d(S,a_j)}{d(S,a_k)}\,. 
\]
Otherwise (note points $A_1$ or $A_2$ in \Figref{fig:denotation}) one of the axes includes the angles $\ssmFrac{\pi + \alpha_i - \gamma_i}2$ or $\ssmFrac{\pi - \alpha_i + \gamma_i}2$ and $\ssmFrac{\pi - \alpha_i-\gamma_i}2$ with $a_j$ and $a_k$, and we obtain
\begin{equation}\label{eq:transmission1}
   \frac{-\sin\gamma_i - \sin\alpha_i}{\sin(\gamma_i + \alpha_i)}
    = -\frac{\cos\ssmFrac{\alpha_i - \gamma_i}2}{\cos\ssmFrac{\alpha_i + \gamma_i}2}
   = \frac{d(S,a_j)}{d(S,a_k)}\,. 
\end{equation}
This confirms the claim that can be rewritten in the form
\begin{equation}\label{eq:transmission2}
  t_1:t_2:t_3 = \frac 1{d(S,a_1)} : \frac 1{d(S,a_2)} : \frac 1{d(S,a_3)}
\end{equation} 
with $\left( d(S,a_1),\,d(S,a_2),\,d(S,a_3)\right)$ as normalized {\em trilinear coordinates} of $S$ w.r.t.\ $A_1A_2A_3$ (see, e.g., \cite[p.~413]{Conics}).
The absolute values of the distances $d(S,a_i)$ are the radii of the concentric circles $k_{12}, k_{23}, k_{31}$ in Figures \ref{fig:Fig_1}a and \ref{fig:denotation}.
\end{Proof}
 
\subsection{Necessary and sufficient condition for flexibility}

The following theorem rephrases a part of Bricard's main result in \cite{Bricard1}.
The new formulation paves the way to higher-dimensional versions. 

\begin{thm}\label{thm2} 
A flat octahedron $A_1\dots\Bq_3$ is twice-flat and flexible in $\Raum E^3$ if and only if it is locally symmetric according to \Defref{def:loc_symm 2d} and at the vertices $A_1$, $A_2$ and $A_3$ one of the axes of symmetry passes through a common (finite or infinite) point $\Sq$.
\end{thm}

\begin{Proof}
(i) We show first that the stated conditions are sufficient for flexibility and begin with the case of a finite $\Sq$:
\\
By virtue of \Lemref{lem:ratio}, every side length of the moving triangle is preserved since for the transmissions between the wings' rotations from $\Wing_1$ to $\Wing_2$, from $\Wing_2$ to $\Wing_3$, and finally from $\Wing_3$ to $\Wing_1$ holds
\begin{equation}\label{eq:prod_ratios}
  \frac{d(\Sq,a_3)}{d(\Sq,a_1)} \cdot \frac{d(\Sq,a_2)}{d(\Sq,a_3)}
   \cdot \frac{d(\Sq,a_1)}{d(\Sq,a_2)}\,t_1 = t_1.
\end{equation}
This confirms the flexibility.

\noindent
If $\Sq$ lies at infinity (see \Figref{fig:parallel_axes}b), then we choose for each $i=1,2,3$ a point $S_i$ on the axis through $A_i$ as the image of $A_i$ under a common translation.
Thus, for each axis $a_k$, $k \ne i$, we obtain equal signed distances $d(S_i,a_k) = d(S_j,a_k)$ if $j\ne k$. 
By \Lemref{lem:ratio} follows again, that the product of the ratios
$d(S_i,a_j)/{d(S_i,a_k)}$ equals 1.

\smallskip\noindent
(ii) Conversely, if a given twice-flat octahedron is flexible, then the transmission ratios between any two wings define by \eqref{eq:transmission2} homogeneous trilinear coordinates $(t_1:t_2:t_3)$ of a unique point $\Sq$ in the plane $\Rasthe$, and by \eqref{eq:transmission1} $\Sq$ is located on one of the axes of local symmetry at each vertex $A_i$ of $\Rast$. 
\end{Proof}

From \eqref{eq:transmission2} follows that at all poses of $\Gang$ the vertices $B_1$, $B_2$ and $B_3$ are corresponding points of projective mappings between the circular paths.
This holds since $t_i = \tan\ssmFrac{\phi_i}2$ defines a projective scale on the trajectory of $B_i$ (see, e.g., \cite[p.~259]{Conics}). 

\begin{cor}\label{cor:concurrent} 
Let $A_1\dots \Bq_3$ be a flat pose of a twice-flat octahedron.
If at the vertices $A_1$, $A_2$ and $A_3$ one of the axes of symmetry passes through a common point $\Sq$, then the same is true for the remaining vertices $\Bq_1$, $\Bq_2$ and $\Bq_3$. 
\end{cor}

\begin{Proof}
Local symmetry and concurrent axes of symmetry at $A_1,A_2,A_3$ imply concentric incircles as shown in \Figref{fig:denotation}.
Note also below the statement (i) in \Thmref{thm:flat_flexible}.
\end{Proof}

\begin{rem}\label{Rem1}
The octahedra with $\Sq$ at infinity are also particular cases of Bricard's line-symmetric Type~1.
Their flexions are determined by the flexions of one four-sided pyramid combined with a line reflection in the axis of symmetry of the bounding isogram. 
These octahedra have the additional property that they admit bifurcations at the flat poses since both axes of symmetry at the vertices belong to families of concurrent axes, which are sufficient for flexiblity by \Thmref{thm2}.
Such bifurcation do not exist in the cases with finite $\Sq$.
\end{rem}

\begin{figure}[htb] 
  \psone{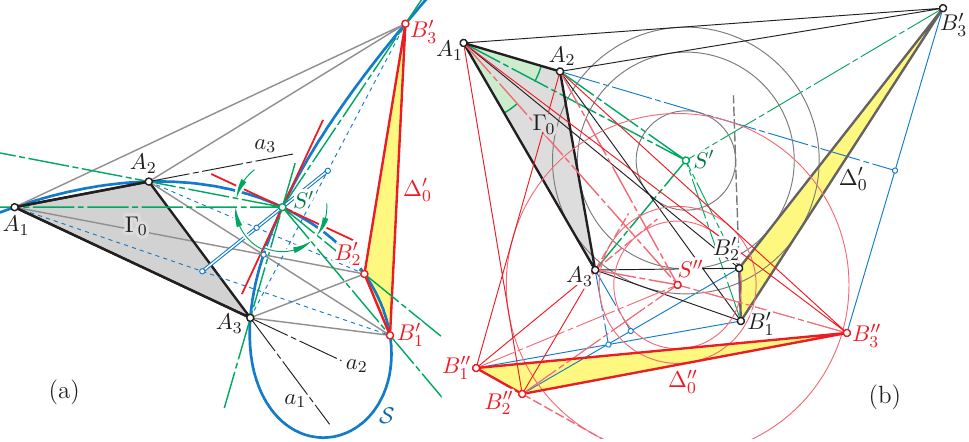}{166mm}
 \caption{The flat poses $A_1\dots\Bq_3$ and $A_1\dots\Bqq_3$ of a flexible twice-flat octahedron:\\
(a) The opposite vertices $(A_i,\Bq_i)$ for $i=1,2,3$ are associated points of a strophoid $\Stroph$ with the node $\Sq$.
The lines $[\Sq,A_i]$, $[\Sq,\Bq_i]$ are symmetric with respect to the tangents of the strophoid at $\Sq$.
The points $[A_i,\Bq_j]\cap[A_j,\Bq_i]$ for $i\ne j$ belong to $\Stroph$, too.
\\
(b) The meeting points $\Sq$ and $\Sqq$ of the concurrent axes of symmetry in the two flat poses are isogonal with respect to the fixed triangle $A_1A_2A_3$.}
\label{fig:strophoid}\label{fig:isogonal}
\end{figure}

\begin{thm}\label{thm:flat_flexible} 
The flat poses of a flexible twice-flat octahedron $\Ccal 3$ have the following properties in addition to those listed above in \Lemref{lem:inscribed_conic}, provided that the point $\Sq$ is finite:
\\[0.5mm]
{\em (i)} The quadrangles $A_1A_2\Bq_1\Bq_2$, $A_2A_3\Bq_2\Bq_3$, and $A_3A_1\Bq_3\Bq_1$ have incircles with the common center $\Sq$ (Figures~\ref{fig:Fig_1}a and \ref{fig:isogonal}b).
The absolute values of the normalized trilinear coordinates of $\Sq$ w.r.t.\ $A_1A_2A_3$ and $B_1B_2B_3$ are equal and proportional to the angular velocities of the wings relative to $\Rast$ at the initial flat pose. 
The same holds for all other pairs of opposite triangles of $\Ccal 3$.
\\[0.5mm]
{\em (ii)} All pairs of opposite vertices $(A_i,\Bq_i)$ are associated points of a strophoid $\Stroph$ with the node $\Sq$ (\Figref{fig:strophoid}a).
\\[0.5mm]
{\em (iii)} The connecting lines $[\Sq,A_i], [\Sq,\Bq_i]$ for $i=1,2,3$ have the tangents to $s$ at the node $\Sq$ as common axes of symmetry.
This results in congruent angles at $\Sq$.
\\[0.5mm]
{\em (iv)} Let $\Sqq$ be the meeting point of the axes of symmetry in the second pose. Then $\Sq$ and $\Sqq$ are isogonal w.r.t.\ the fixed triangle $A_1A_2A_3$ (\Figref{fig:isogonal}b).
\end{thm}

\begin{Proof}
(i) The three quadrangles have angle bisectors meeting at $\Sq$.
The proportion \eqref{eq:transmission2} is also valid for the initial angular velocities $\dot t_1(0):\dot t_2(0):\dot t_3(0)$ due to $\dot\phi_i(0) = 2\,\dot t_i(0)$. 

\smallskip\noindent
(ii) We recall that the strophoid $\Stroph$ is the geometric locus of points $X$ with the property that one angle bisector between the connections with $A_1$ and $\Bq_1$ passes through the node $\Sq$.
The tangents at the node bisect the connections with $A_1$ and $\Bq_1$, which are the fixed lines of the symmetric involution. 
Any two lines through $\Sq$ and symmetric w.r.t.\ the node tangents intersect $s$ additionally at a pair of associated points (note \cite{Sta_156}).
As a consequence, the point of intersection between $[A_i,\Bq_j]$ and $[A_j,\Bq_i]$ for $i\ne j$ is also located on $\Stroph$.

\smallskip\noindent
(iii) For each quadrangle, the pairs of line pencils with carriers at opposite vertices and the incircle belong to a dual pencil of conics.
The respective tangents through $\Sq$ define a symmetric involution.
Any two quadrangles share a pair of opposite vertices.
Therefore, the three involutions are equal.

\smallskip\noindent
(iv) We refer to (\Figref{fig:isogonal}b).
Since $[A_1,\Bq_3]$ and $[A_1,\Bqq_3]$ are symmetric w.r.t.\ $[A_1,A_2]$, the respective angle bisectors with $[A_1,A_3]$, i.e., the lines $[A_1,\Sq]$ and $[A_1,\Sqq]$ include with $[A_1,A_3]$ the angles $(\alpha_1 + \gamma_1)/2$ and $(\alpha_1 - \gamma_1)/2$ according to the notation used in \Figref{fig:denotation}. 
When replacing $\Sq$ by $\Sqq$, the ratio in \Lemref{lem:ratio} becomes reciprocal, i.e., $d(\Sq,a_2):d(\Sq,a_3) = d(\Sqq,a_3):d(\Sq,a_2)$, since the angles $\phi_i$ of rotations are replaced by $\phi_i - \pi$. 
\end{Proof}

\begin{rem}\label{Rem2}
At the flat pose $A_1\dots\Bq_3$ of a twice-flat octahedron, the three pairs of line pencils with carriers $(A_i, \Bq_i)$ span a two-parameter linear set (or net) of dual conics.
This net contains the conics mentioned in \Lemref{lem:inscribed_conic},\,(i).
If the octahedron is flexible, then by virtue of \Thmref{thm:flat_flexible},\,(i) the net contains also the concentric circles centered at $\Sq$, provided that $\Sq$ is finite.
Then, the net contains the two-fold counted line pencil with the carrier $\Sq$ as well as the set of isotropic lines and, consequently, with each conic all confocal conics.
The strophoid $\Stroph$ is Cayley's curve of this net (see, e.g., \cite[p.~334]{Conics}).
Herewith, we obtain a characterization of the configuration of vertices at flat poses of type-3 octahedra. 
\end{rem}

\medskip
How to construct a flat pose of a flexible twice-flat octahedron?
Bricards first construction in \cite[p.~144]{Bricard1} is based on \Thmref{thm2} and slightly different from that depicted in \Figref{fig:Fig_1}a:
It begins with the choice of a triangle $A_1A_2A_3$ and a point $\Sq$ outside all extended sides $[A_i,A_j]$.
Then the line $[A_i,B_j]$, $i\ne j$, is obtained by reflection of $[A_i,A_j]$ in $[A_i,\Sq]$.
Thus, point $B_j$ is determined as the intersection of two lines.
It will turn out that a similar construction is valid for flexible twice-flat cross-polytopes in $\Raum E^n$.  

\subsection{Spatial poses of flexible twice-flat octahedra}

The coming theorem enumerates for flexible twice-flat octahedra some properties of their spatial poses, in particular local symmetries.
These symmetries do not act on faces but on the planes spanned by faces.
It is noteworthy that many of these properties can already be extracted from Figures~5 and 6 in Bennett's paper \cite{Bennett}.
Below we use for all permutations $(i,j,k)$ of $(1,2,3)$ beside $\Wing_i = B_iA_jA_k$ the symbol $\CoWing_i:= A_iB_jB_k$ for the wings of the inverse motion $\Rast/\Gang$.

\begin{thm}\label{thm:spatial_pose} 
Let $A_1\dots B_3$ be a spatial pose of a flexible twice-flat octahedron $\Ccal 3$ with finite centers $\Sq, \Sqq$ in the planar poses.
Then holds for each permutation $(i,j,k)$ of $(1,2,3)$:
\\[0.5mm]
{\em (i)} The dihedral angle $\psi$ between the faces $X_iX_jA_k$ and $X_iX_jB_k$ through the edge $X_iX_j$ with $X_i\in\{A_i,B_i\}$ and $X_j\in\{A_j,B_j\}$ is either congruent or supplementary to the angle $\phi_k$ between the wing $\Wing_k$ and the fixed triangle $\Rast$.
Congruence holds if and only if at the initial flat pose the dihedral angle $\psi$ equals $0$.
\\[0.5mm]
{\em (ii)} The pairs of lines $([A_iA_j],\, [A_iB_j])$ and $([A_iA_k],\, [A_iB_k])$ through $A_i$ have a common bisecting plane $\mu_i$ which is also common to the pairs of planes $([A_iA_jA_k],\,[A_iB_jB_k])$ and $([A_iA_jB_k],\, [A_iB_jA_k])$.
The three bisecting planes $\mu_1,\mu_2,\mu_3$ share a line $s$.
The analogous planes of symmetry $\nu_i$ through $B_i$ pass through $s$, too. 
\\[0.5mm]
{\em (iii)} The sides of the three skew quadrangles $A_1A_2B_1B_2$, $A_2A_3B_2B_3$ and $A_3A_1B_3B_1$ are located respectively on one-sheeted hyperboloids of revolution $\mathcal H_1, \mathcal H_2, \mathcal H_3$ with the axis $s$ (\Figref{fig:hyps}).
The three diagonals $[A_1,B_1], \dots,[A_3,B_3]$ belong to a hyperbolic paraboloid.
\\[0.5mm]
{\em (iv)} For all $i$, the pairs of planes $(\mu_i,\nu_i)$ connecting $s$ with $A_i$ and $B_i$ have common planes of symmetry. 
\\[0.5mm]
{\em (v)} The midpoints of the diagonals $A_iB_i$ are coplanar with $s$, and the center of each hyperboloid $\mathcal H_i$ is aligned with two of these midpoints.  
\\[0.5mm]
{\em (vi)} The vertices $A_1,\dots,B_3$ belong to a spatial circular cubic $\Stroph$ with $s$ as a bisecant such that the tangents to $\Stroph$ at the meeting points with $s$ lie in orthogonal planes through $s$. 
\\[0.5mm]
{\em (vii)} The four planes spanned by $\Gang$, $\Wing_1$, $\Wing_2$, and $\Wing_3$ are tangent to a cone of revolution with the axis $s$.
The same holds for the planes spanned by $\Rast,\dots,\CoWing_{3}$ (\Figref{fig:hyps}).
The rotation about $s$ with $[\Wing_i]\mapsto[\Wing_j]$ takes also $[\CoWing_j]$ to $[\CoWing_i]$ for all $i,j\in\{0,\dots3\}$.
The two quadruples of planes form a pair of Moebius tetrahedra.
\end{thm}

\begin{figure}[htb] 
  \psone{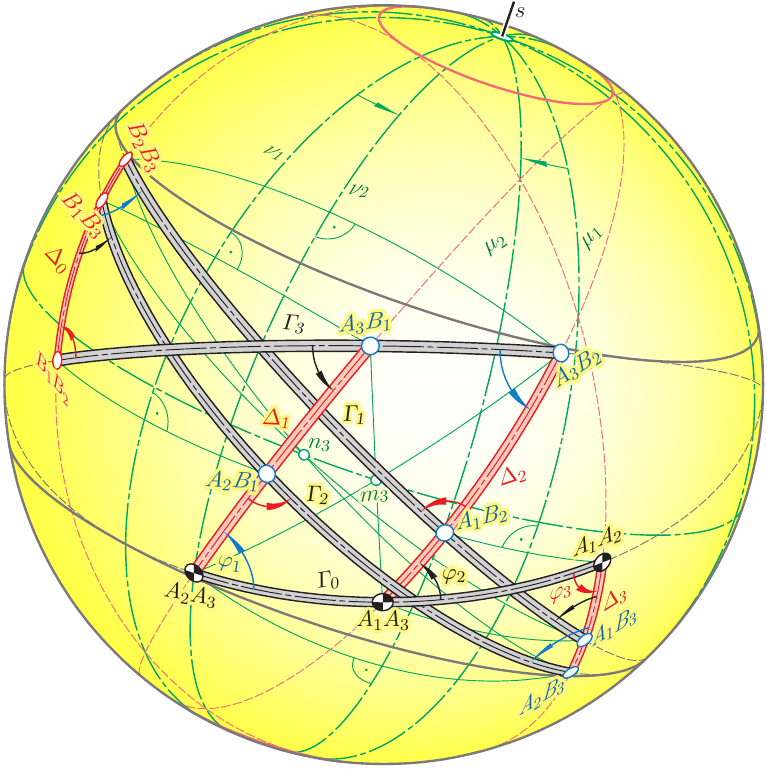}{130mm}
  \caption{The spherical image of the faces and edges of a type-3 octahedron is an 8-bar linkage with 12 joints \cite{Bennett}.
 The planes of the wings and $[\Gang]$ are equally inclined w.r.t.\ an axis $s$.
The same holds for the planes spanned by the faces $\CoWing_0\,\dots,\CoWing_3$.}
  \label{fig:spher_linkage2}
\end{figure}

\begin{Proof}
(i) This holds since at each isogonal double-pyramid with apex out from $\{A_1, \dots, B_3\}$ opposite dihedral angles are congruent (note the octahedron's spherical image in \Figref{fig:spher_linkage2}).
It implies, that along each quadrangle $A_jA_kB_jB_k$ the dihedral angles are congruent or supplementary.

\smallskip\noindent
(ii) By virtue of \Lemref{lem:symm_doublePyr}, the isogonal double-pyramid with apex $A_i$ that extends the pyramid through the quadrangle $A_jA_kB_jB_k$ has a plane $\mu_i$ of symmetry.
The reflection 
in $\mu_i$ preserves the pyramid's diagonal planes $[A_i,A_j,B_j]$, $[A_i,A_k,B_k]$ and exchanges within these planes the lines $[A_iA_j]$ and $[A_iB_j]$ as well as $[A_iA_k]$ and $[A_iB_k]$ (note \Figref{fig:double-cone},\,right, and choose $A_i=O$, $A_j\in g_1$, $A_k\in g_2$, $B_j\in g_3$, and $B_k\in g_4$). 
Moreover, the reflection in $\mu_i$ exchanges opposite planes, i.e., 
\[ [\Wing_k] = [A_iA_jB_k]\,\leftrightarrow\,[A_iB_jA_k] = [\Wing_j], \quad   
    [\Rast] = [A_iA_jA_k] \,\leftrightarrow\,[A_iB_jB_k] = [\CoWing_i].
\] 
The plane $\mu_i$ connects the intersection lines $[\Wing_j]\cap[\Wing_k]$ and $[\Rast]\cap[\CoWing_i]$.
\\
Similarly, the opposite double-pyramid with apex $B_i$ has a plane $\nu_i$ of symmetry.
The reflection in $\nu_i$ exchanges $[B_iA_j]$ and $[B_iB_j]$ as well as $[B_iA_k]$ and $[B_iB_k]$, and moreover 
\[ [\CoWing_j] = [B_iA_jB_k]\,\leftrightarrow\,[B_iB_jA_k] = [\CoWing_k], \quad  
    [\Wing_i] = [B_iA_jA_k] \,\leftrightarrow\,[B_iB_jB_k]=[\Gang].
\]
Let $s$ be the line of intersection $\mu_1\cap\mu_2$ of the planes of symmetry through $A_1$ and $A_2$.
The line $s$ must be finite since otherwise the respective mirrors $[\Wing_1]$ and $[\Wing_2]$ of $[\Wing_3]$ would be parallel which contradicts (i).  
\\ 
Each point $S\in s$ has equal distances to $[\Rast]$ and $[\CoWing_1]$ as well as to $[\Rast]$ and $[\CoWing_2]$. 
The same holds for the distances to $[\Wing_3]$ and $[\Wing_2]$ as well as to $[\Wing_3]$ and $[\Wing_1]$.
\\
From equal distances to $[\CoWing_1]=[A_1B_2B_3]$ and $[\CoWing_2]=[B_1A_2B_3]$, which are opposite planes through $B_3$, follows that $s$ belongs either to $\nu_3$ or to the second bisecting plane of these two planes.
The latter can be excluded since in the initial flat pose $s$ coincides with the perpendicular to $\Rasthe$ through $\Sq$ or $\Sqq$ and all $\mu_i$ and $\nu_i$ pass through $s$, and continuity guarantees that the bisecting planes cannot switch during the flexion. 
Thus, the plane $\nu_3$ passes through $s$, too.
Similarly, equal distances to the opposite planes $[\Wing_1]=[B_1A_2A_3]$ and $[\Wing_2]=[A_1B_2A_3]$ through $A_3$ imply $s\subset\mu_3$. 
Iteration by cyclically increasing all involved subscripts confirms the statement.

\smallskip\noindent
(iii) Concerning the skew quadrangle $A_1A_2B_1B_2$, the consecutive reflections in the planes $\mu_2$, $\nu_1$, $\nu_2$, and $\mu_1$ take $[A_1,A_2]$ via $[A_2,B_1]$ and  $[B_1,B_2]$ to $[B_2,A_1]$. 
This confirms the existence of the hyperboloid $\mathcal H_3$ passing through the four lines.
It agrees with the fact that due to the existing incircles in the flat poses, an alternating sum of the side lengths vanishes, which guarantees the location of the skew quadrangle on a hyperboloid of revolution (note \cite[Theorem~2.2.2]{Quadrics}).
\\
As noted in (ii), the diagonal planes of the double-pyramid with apex $A_i$ are orthogonal to $\mu_i$.
Hence, also the line of intersection $m_i:= [A_i,A_j,B_j] \cap [A_i,A_k,B_k]$ through $A_i$ is orthogonal to $\mu_i$ and, consequently, also orthogonal to $s$.
Note that $m_i$ is the transversal drawn from $A_i$ to the remaining diagonals $[A_j,B_j]$ and $[A_k,B_k]$.  
The lines $n_i\perp\nu_i$ through $B_i$ play a similar role.
Common transversals of the three diagonals $[A_i,B_i]$ for $i=1,2,3$ like $m_i$ and $n_i$ belong to the complementary regulus of the quadric defined by the diagonals.
Since $m_1,\dots, n_3$ are orthogonal to $s$, the quadric through the three diagonals is a hyperbolic paraboloid.

\begin{figure}[htb] 
  \psone{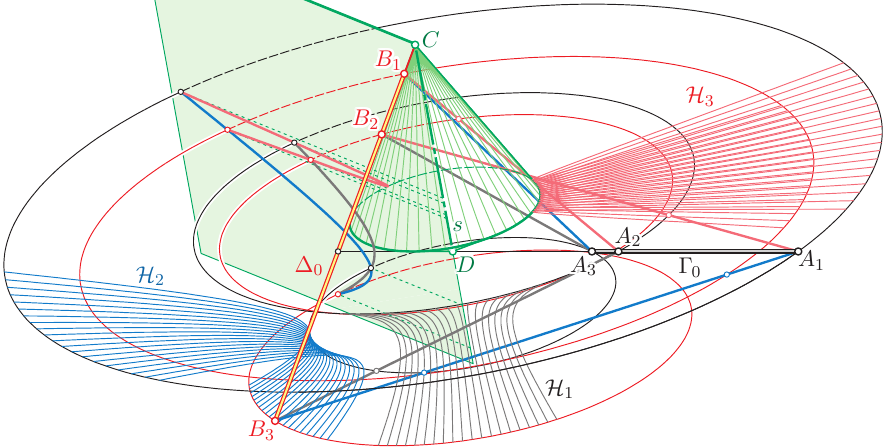}{150mm} 
  \caption{Spatial pose of a flexible twice-flat octahedron showing both triangles $\Rast = A_1A_2A_3$ and $\Gang = B_1B_2B_3$ in an edge view.
The three quadrangles $A_2A_3B_2B_3$ (grey), $A_3A_1B_3B_1$ (blue) and $A_1A_2B_1B_2$ (red) are located on hyperboloids of revolution $\mathcal H_1, \dots, \mathcal H_3$ with the common axis $s$ and displayed meridians. 
  The planes $[\Rast], \dots, [\CoWing_3]$ envelop a cone of revolution with axis $s$ and apex $D\in [\Rast]$.
Similarly, $[\Gang]$ and the planes spanned by the wings envelop a cone with axis $s$ and apex $C\in [\Gang]$.}
  \label{fig:hyps} 
\end{figure}

\smallskip\noindent
(iv) The quadrics through the sides of the skew quadrangle $A_iA_jB_iB_j$ belong to a dual pencil (range).
This range includes the pairs of bundels of planes with carriers $(A_i,B_i)$ and $(A_j,B_j)$.
According to Desargues's involution theorem, each dual quadric of the range sends through $s$ two planes that define an involution.
Since the hyperboloid of revolution sends two isotropic tangent planes through $s$, the involution is symmetric.
This symmetric involution is the same for all three quadrangles, as any two share a pair of opposite vertices.

\smallskip\noindent
(v) The quadrics of the range mentioned before have their centers on a line passing through the midpoints of the diagonals $A_iB_i$ and $A_jB_j$. 
This center line intersects $s$ at the center of hyperboloid of revolution $\mathcal H_k$ through the quadrangle.
Since any two of the three quadrangles share a pair of opposite vertices, the center lines of the three ranges must be coplanar with the axis $s$.

\smallskip\noindent
(vi) The cubic is defined as locus of points $X$ with the property that the connecting lines with two given points $A_1,B_1$ have a plane of symmetry passing through $s$ (see \cite{Bricard2}).\footnote{
Points $X\in\Stroph$ can be found as follows:
For any given plane $\eps$ through $s$, reflect $A_1$ in $\eps$, connect the image with $B_1$ and intersect the line with $\eps$. 
The mirror points are located on a circle through $A_1$ with axis $s$.}

\smallskip\noindent
(vii) Referring to (ii), the reflection $\sigma_i$ in the plane $\mu_i$ through $s$ exchanges $[\Rast]$ with $[\CoWing_{i}]$ for $i=1,2,3$.
Consequently, the planes $[\Rast],\dots,[\CoWing_3]$ are equally inclined w.r.t.\ $s$ and share a point $D\in s$.
Thus, they envelop a cone of revolution with axis $s$ and apex $D$.
Analogously, the planes $[\Gang],\dots,[\Wing_3]$ envelope a second cone with axis $s$ and apex $C$ (\Figref{fig:spher_linkage2}).
Moreover, the planes $[\Rast],[\Wing_1],[\Wing_2],[\Wing_3]$ enclose the tetrahedron $A_1A_2A_3C$ that is in- and circumscribed the tetrahedron $B_1B_2B_3D$ bounded by $[\Gang],[\CoWing_1],[\CoWing_2],[\CoWing_3]$. 
\\
By the way, the infinitesimal screw motion of $\Gang$ relative to $\Rast$ is combined with a null polarity that maps $B_i$ to the orthogonal plane $[\Wing_i]$ of its trajectory, while $C$ is the null point of $\Gang$. 
This null polarity interchanges the two tetrahedra.
\end{Proof}

\begin{rem}\label{rem:Grundriss}
The properties listed in \Thmref{thm:spatial_pose} reveal, that the orthogonal projection of a spatial pose in direction of the axis $s$ yields an image that has all properties of the flat poses as listed in \Thmref{thm:flat_flexible}.
The spatial cubic $\Stroph$ from (vi) is sent to the strophoid from \Thmref{thm:flat_flexible},\,(ii).
The negative pedal curve of the strophoid w.r.t.\ its node is a parabola (see \cite[Fig.~2.27]{Conics}) that equals the visual contour of the paraboloid through the diagonals $[A_i,B_i]$ and the lines $m_1,\dots,n_3$. 
\end{rem}

At the spherical image of a spatial pose the octahedron's faces and edges are represented as bars and joints of a spherical 8-bar linkage with 12 revolute joints as the union of six isogonal spherical four-bar linkages \cite{Bennett}.
Three four-bars have their bases on the image of $\Rast$, while their respective couplers are aligned with arms of other three four-bars with bases on the spherical image of $\Gang = B_1B_2B_3$ (\Figref{fig:spher_linkage2}).

\section{Flexible twice-flat cross-polytopes in $\Raum E^n$}

As reported in Section~\ref{sec:Intro}, a cross-polytope $\Ccal n$ in $\Raum E^n$ has $n$ pairs of opposite vertices $(A_i,B_i)$, $i=1,\dots,n$.
We return to the original denotation and assume that the simplex $\Rastn:= A_1\dots A_n$ in the hyperplane $\Rasthe$ is fixed while the opposite simplex $\Gangn:= B_1\dots B_n$ in the hyperplane $\Ganghe$ is moving.
Hence, during the rotation of the simplex called wing $\Wing_i:= A_1 \dots A_{i-1} B_i A_{i+1} \dots A_n$ about the $(n\minus 2)$-dimensional hinge $\Achse_i = [A_1 \dots A_{i-1} A_{i+1} \dots A_n]$ the vertex $B_i$ traces a circle (note case $n=3$ in Figs.~\ref{fig:in_motion} and \ref{fig:moving} and case $n=4$ in \Figref{fig:4D-cross-polytope}).

Two facets of $\Ccal n$ like $X_1\dots A_i\dots X_n$ and $X_1\dots B_i\dots X_n$ with $X_k\in\{A_k,B_k$ are called {\em neighboring}, if they share an $(n\minus 2)$-face.
In the case of a flexible cross-polytope, the relative motion of two neighboring facets is a rotation about the hinge spanned by the common $(n\minus 2)$-face.
It will turn out that the $(n\minus 3)$-faces of $\Ccal n$ serve as apices of the higher-dimensional versions of the isogonal pyramids at the type-3 octahedra.
Each $(n\minus 3)$-face of $\Ccal n$ is the meet of four facets where each two consecutive facets in cyclic order are neighboring. 
For example, the facets through the $(n\minus 3)$-face $X_3\dots X_n$ are
$A_1A_2X_3\dots X_n$, $A_1B_2X_3\dots X_n$, $B_1B_2X_3\dots X_n$, and $B_1A_2X_3\dots X_n$.  

\subsection{Local symmetries at flat poses of $\Ccal n$}

Let $\Ccal n$ be a twice-flat cross-polytope in $\Raum E^n$ for $n>3$.
This means by \Defref{def:twice-flat} that $\Ccal n$ admits beside a given flat pose $A_1\dots A_n \Bq_1 \dots \Bq_n$ a second flat pose $A_1\dots A_n\Bqq_1\dots\Bqq_n$ with $\Bq_i\ne \Bqq_i$ for all $i\in\{1,\dots,n\}$.

\begin{figure}[htb] 
  \psone{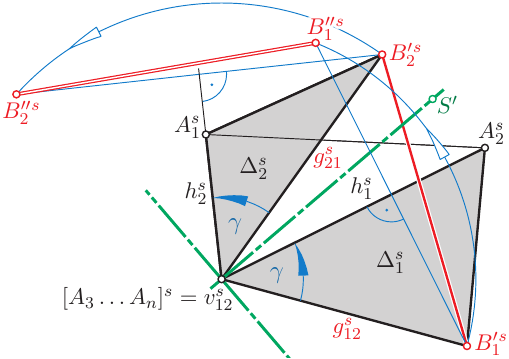}{85mm}  
  \caption{Two-dimensional view of the two wings $\Wing_1, \Wing_2$ in a flat pose within $\Rasthe$ after an orthogonal projection of $\Rasthe$ parallel to $\Kante_{12} = [A_3\dots A_n]$ into a 2-plane.}
  \label{fig:proj}
\end{figure}

The second pose $\Bqq_1$ of the vertex $B_1$ arises in $\Raum E^n$ from the first pose $\Bq_1$ by a halfturn about $\Achse_1$.
This halfturn acts within the fixed hyperplane $\Rasthe$ like a reflection in $\Achse_1$.
Consequently, within $\Rasthe$ the $(n\minus 2)$-dimensional space $\Achse_1$ is the bisecting hyperplane between $\Bq_1$ and $\Bqq_1$. 
Similarly, $\Achse_2$ is the bisecting hyperplane between $\Bq_2$ and $\Bqq_2$.
Because of equal distances $\ol{\Bq_1\Bq_2} = \ol{\Bqq_1\Bqq_2}$ we obtain congruent simplices $\Bq_1\Bq_2A_3\dots A_n$ and $\Bqq_1\Bqq_2A_3\dots A_n$.
Hence, there exists in $\Rasthe$ a rotation about the $(n\minus 3)$-dimensional space $\Kante_{12}:= \Achse_1\cap\Achse_2 = [A_3 \dots A_n]$ with $\Bq_1\mapsto\Bqq_1$ and $\Bq_2\mapsto\Bqq_2$ (\Figref{fig:proj}). 

Within $\Rasthe$, the rotation about $\Kante_{12}$ through the halved angle takes $[A_3 \dots A_n \Bq_1] = [\Kante_{12},\Bq_1]$ to the bisecting hyperplane $\Achse_1 = [\Kante_{12},A_2]$ and $[\Kante_{12},\Bq_2]$ to $\Achse_2 = [\Kante_{12},A_1]$.
Thus, within $\Rasthe$ the pairs of hyperplanes $([\Kante_{12},A_1], [\Kante_{12}, \Bq_1])$ and $([\Kante_{12},A_2], [\Kante_{12},\Bq_2])$ have a common bisecting hyperplane.

Similar symmetries exist at all other pairs of vertices $(A_i,\Bq_i)$ and $(A_j,\Bq_j)$, $i\ne j$, for their $(n\minus 2)$-dimensional connections with the $(n\minus 3)$-dimensional space $\Kante_{ij}$ spanned by the remaining vertices $A_k$ of $\Rastn$ with $k\ne i,j$. 

\begin{thm}\label{thm:twice-flat n} 
A flat cross-polytope $\Ccal n$ with the simplex $\Rastn = A_1\dots A_n$ and the coplanar opposite simplex $\Gangnq = \Bq_1\dots \Bq_n$ is twice-flat if and only if for all $i,j\in\{1,\dots,n\}$ the connecting $(n\minus 2)$-spaces of the pairs $(A_i,\Bq_i)$ and $(A_j,\Bq_j)$ with the $(n\minus 3)$-face of $\Rastn$ opposite to $A_iA_j$ have common angle bisecting $(n\minus 2)$-spaces.
\end{thm}

\begin{Proof}
Let from now on $\Kante_{ij}=\Kante_{ji}$ denote the span of $\left\{A_k\,|\,k\in\{1\dots n\}\setminus\{i,j\}\right\}$.
Then, the stated condition for being twice-flat is necessary and sufficient in $\Rasthe$ as common angle bisectors of the pairs of hyperplanes $([\Kante_{ij}, A_i], [\Kante_{ij}, \Bq_i])$ and $([\Kante_{ij}, A_j], [\Kante_{ij}, \Bq_j])$ characterize equal distances $\ol{\Bq_i\Bq_j} = \ol{\Bqq_i\Bqq_j}$.
\end{Proof} 

Since for the existence of a second flat pose it is not relevant which facet of the cross-polytope $\Ccal n$ is fixed, the symmetries stated above result in similar symmetries at all $(n\minus 3)$-faces of $\Ccal n$.

\begin{deft}\label{def:loc_symm dn} 
A flat cross-polytope $A_1\dots \Bq_n$ in a hyperplane of $\Raum E^n$ is called {\em locally symmetric} when it satisfies the condition mentioned in \Thmref{thm:twice-flat n}.
In this case similar symmetries hold for each $(n\minus 3)$-face of the flat cross-polytope.
\end{deft}

\subsection{Transmission between two wings of a twice-flat cross-polytope}

Now we investigate the relation between the angles of rotation $\phi_1, \phi_2$ of the wings $\Wing_1 = A_2A_3\dots A_nB_1$ and $\Wing_2 = A_1A_3\dots A_nB_2$ relative to the fixed simplex $\Rastn = A_1\dots A_n$ when the moving edge $B_1B_2$ preserves its length.

The wing $\Wing_1$ with $B_1$ rotates about $\Achse_1 = [A_2A_3\dots A_n]$, while $\Wing_2$ with $B_2$ rotates about $\Achse_2 = [A_1A_3\dots A_n]$.
The two wings are connected by the facet $\Wing_{12}:= A_3\dots A_n B_1 B_2$.
In non-flat poses, the four facets $\Wing_1$, $\Wing_{12}$, $\Wing_2$, and $\Rastn$ form a four-sided pyramid in $\Raum E^n$ with the $(n\minus 3)$-dimensional apex space $\Kante_{12} = [A_3\dots A_n] = \Achse_1\cap\Achse_2$ and the quadrangular base $B_1B_2A_1A_2$.
The `edges' of the pyramid are the intersections between adjacent facets, i.e., the $(n\minus 2)$-faces $\Wing_1\cap\Wing_{12} = A_3\dots A_n B_1$, $\Wing_{12}\cap\Wing_2 = A_3\dots A_n B_2$, $A_3\dots A_n A_1$ spanning $\Achse_2$, and $A_3\dots A_n A_2$ spanning $\Achse_1$.

\smallskip
For the investigation of the pyramid's flexion in $\Raum E^n$, we recall from Linear Algebra that two (linear) subspaces $U, V$ in $\Raum E^n$ are called {\em orthogonal} (or total-orthogonal) if each direction in $U$ is orthogonal to all directions in $V$.
If the sum of dimensions of the two spaces equals $n$, then conversely each direction that is orthogonal to all directions in $V$ is parallel to $U$.

\begin{lem}\label{lem:total-orthogonal} 
Let $(U_i,V_i)$ for $i=1,2$ be two pairs of orthogonal subspaces in $\Raum E^n$ where in both cases the dimensions sum up to $n$.
Then $U_1\subset U_2$ implies that each direction in $V_2$ is parallel to $V_1$.
\end{lem}

\begin{figure}[htb] 
  \psone{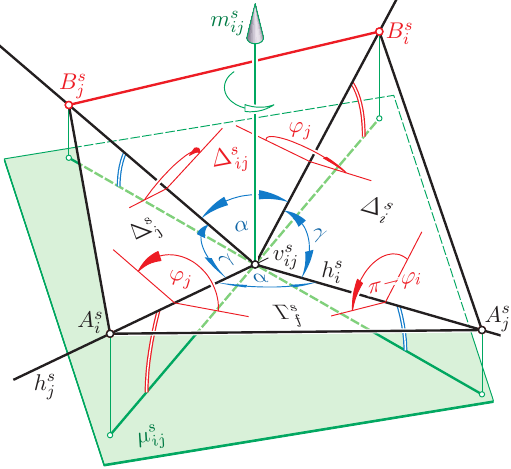}{85mm} 
  \caption{The orthogonal projection $p_{ij}$ of $\Raum E^n$ parallel to $\Kante_{ij}=\Achse_i\cap\Achse_j$ into a 3-space $\Pi_{ij}$ yields this view of the isogonal pyramid with the $(n\minus 3)$-dimensional apex space $\Kante_{ij}$ and the two wings $\Wing_i, \Wing_j$ in a spatial pose.
The plane $\mu_{ij}^s$ is the image of a hyperplane $\mu_{ij}$ of symmetry of the extended double-pyramid.
The preimage of $m_{ij}^s$ is the $(n\minus 2)$-dimensional axis $m_{ij}$ of a halfturn that exchanges opposite elements of the double-pyramid in $\Raum E^n$, where $m_{ij}$ and $\mu_{ij}$ are orthogonal spaces in the bundle with carrier $\Kante_{ij}$.}
  \label{fig:proj_b}
\end{figure}

\begin{figure}[htb] 
  \psone{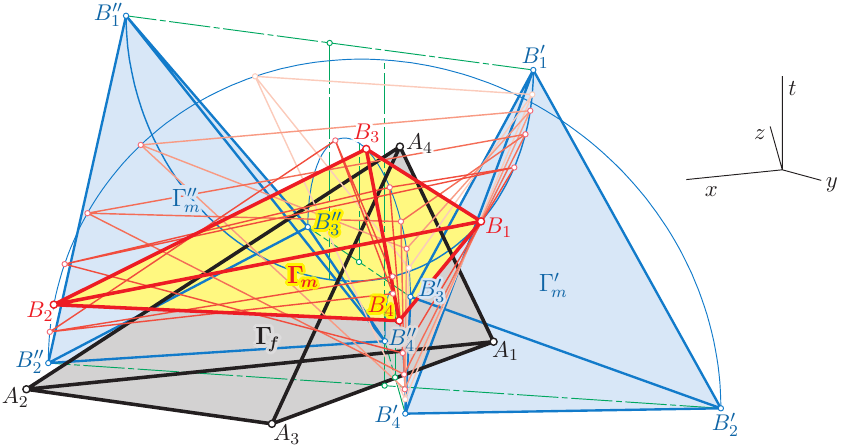}{140mm} 
  \caption{Two-dimensional orthogonal view of a four-dimensional flexible twice-flat cross-polytope $\Ccal 4$ with different poses of the moving tetrahedron $\Gangn = B_1\dots B_4$ in $\Raum E^4$ relative to the fixed tetrahedron $\Rastn = A_1\dots A_4$ in the hyperplane $t = 0$.}
\label{fig:4D-cross-polytope}
\end{figure}

Our analysis utilizes a projection $p_{12}$ of $\Raum E^n$ in direction of the $(n\minus 3)$-dimensional apex space $\Kante_{12}$ into a 3-space $\Pi_{12}$, which is orthogonal to $\Kante_{12}$.
We call the image of this orthogonal projection briefly a {\em side view} and mark images with the superscript $s$.
Thus, the side view of the pyramid in $\Raum E^n$ is a four-sided pyramid in $\Pi_{12}$ with the point $\Kante_{12}^s$ as apex and the quadrangular base $B_1^s B_2^s A_1^s A_2^s$ (\Figref{fig:proj_b}).
The restriction of $p_{12}$ to the hyperplane $\Rasthe = [\Rastn]$ sends the pyramid's flat pose in a 2-plane $\Rasthe^s\subset\Pi_{12}$ (\Figref{fig:proj}).  

\begin{lem}\label{lem:projection} 
The orthogonal projection $p_{12}\!: \Raum E^n \to \Pi_{12}$ in direction of the $(n\minus 3)$-dimensional space $\Kante_{12} = [A_3\dots A_n]$ maps the four-sided pyramid through $B_1B_2A_1A_2$ with the apex space $\Kante_{12}$ to a three-dimensional pyramid while dihedral angles between facets and interior angles at the apex $\Kante_{12}^s$ are preserved. 
The projection $p_{12}$ preserves also the distances of points $X\in\Rasthe$ to the hinges $\Achse_1$ and $\Achse_2$.
\end{lem}

\begin{Proof}  
(i) The dihedral angle $\phi_1$ between the wing $\Wing_1$ and the fixed facet $\Rastn$ is measured in a 2-plane orthogonal to the extended edge $[\Wing_1\cap\Rastn] = \Achse_1$.
From $\Kante_{12}\subset \Achse_1$ follows by \Lemref{lem:total-orthogonal} that the 2-plane with the angle $\phi_1$ lies parallel to the image space $\Pi_{12}$. 
Consequently, the side view shows the dihedral angle $\phi_1$ in true size.
The same holds for the other dihedral angles of the pyramid in $\Raum E^n$ (\Figref{fig:proj_b}).
\\[0.5mm]
(ii) Within the wing $\Wing_1$, the dihedral angle between the two extended edges $[\Kante_{12},A_2] = \Achse_1$ and $[\Kante_{12},B_1]$ is measured in a 2-plane that is orthogonal to $\Kante_{12}$ and therefore parallel to $\Pi_{12}$.
Since the same holds for the other three facets of the four-sided pyramid in $\Raum E^n$, the orthogonal projection $p_{12}$ preserves the pyramid's interior angles at the apex (\Figref{fig:proj_b}).
\\[0.5mm]
(iii) Within $\Rasthe$, the distance of any point $X$ to the $(n\minus 2)$-dimensional hinge $\Achse_1$ is measured along a line that is orthogonal to $\Achse_1$.
From $\Kante_{12}\subset\Achse_1$ follows within $\Rasthe$ by \Lemref{lem:total-orthogonal} that this line is parallel to $\Rasthe^s\subset\Pi_{12}$.
Therefore, the side view shows for each point $X\in\Rasthe$ the distance to $\Achse_1$ in true size as well that to $\Achse_2$.
\end{Proof}

\begin{rem}\label{rem:cylinder}
\Lemref{lem:projection} reveals that the four-sided pyramid in $\Raum E^n$ with the apex space $\Kante_{12}$ can also be seen as a prism with generators parallel to $\Kante_{12}$ and the three-dimensional pyramid in $\Pi_{12}$ as an orthogonal cross-section.
\end{rem}

\smallskip
From \Lemref{lem:projection} follows that the transmission from wing $\Wing_1$ to $\Wing_2$ in $\Raum E^n$ satisfies the same formula like that in 3-space.
This is a result that can already be found in \cite[Proposition~3.3]{Gaifullin}.

By virtue of \Thmref{thm:twice-flat n}, the four-sided pyramid in $\Pi_{12}$ is isogonal, i.e., opposite angles at the apex $\Kante_{12}^s$ are congruent or supplementary (\Figref{fig:proj}).
The same holds, due to \Lemref{lem:projection}, for the preimage in $\Raum E^n$, and the relation between the revolute angles $\phi_1$ and $\phi_2$ matches \eqref{eq:transmission0}, where $\alpha$ and $\gamma$ are interior angles at the $n$-dimensional pyramid with apex space $\Kante_{12}$.
In addition, due to \Lemref{lem:projection}, also eq.\ \eqref{eq:transmission0} is still valid for the distances from $\Achse_1$ and $\Achse_2$ within $\Rasthe$.

Similar orthogonal projections for all other $(n\minus 3)$-faces of the fixed simplex $\Rastn$ yield analogous results.
Thus we can state:

\begin{lem}\label{lem:4d-transmission} 
At a twice-flat cross-polytope $\Ccal n$ in $\Raum E^n$ the distance between two vertices $B_i, B_j\in\Gangn$ remains constant if and only if the angles of rotation $\phi_i, \phi_j$ of the wings $\Wing_i$ and $\Wing_j$ satisfy 
\begin{equation}\label{eq:tan_proportional_n}
   t_i:t_j = \tan\frac{\phi_i}2: \tan\frac{\phi_j}2 =
   \frac 1{d(\Sq,\Achse_i)} : \frac 1{d(\Sq,\Achse_j)},
\end{equation}
provided that, in a flat pose $A_1,\dots,\Bq_n$ of $\Ccal n$, the point $\Sq\in \Rasthe$ lies in a common bisecting $(n\minus 2)$-space through $\Achse_i\cap\Achse_j = \Kante_{ij}$ between the $(n\minus 2)$-dimensional spaces $[\Kante_{ij},A_i]$ and $[\Kante_{ij}, \Bq_i]$ as well as between $[\Kante_{ij},A_j]$ and $[\Kante_{ij}, \Bq_j]$ according to \Thmref{thm:twice-flat n}.
\\
The symbol $d(\Sq,\Achse_k)$ with $k = 1,\dots,n$ stands for the signed distance of the point $\Sq\in\Rasthe$ to the hinge $\Achse_k$ of the wing $\Wing_k$, such that for each hinge in $\Rasthe$ the distance is positive if and only if $\Sq$ belongs to the halfspace containing the interior of the fixed simplex $\Rastn = A_1\dots A_n$.
\end{lem}

\subsection{Necessary and sufficient condition for flexibility of $\Ccal n$}

\begin{thm}\label{thm:fund_thm_En} 
A flat cross-polytope $A_1\dots A_n\Bq_1\dots \Bq_n$ in the hyperplane $\Rasthe$ is twice-flat and flexible in $\Raum E^n$, $n > 3$, if and only if it is locally symmetric according to \Defref{def:loc_symm dn} and at every $(n\minus 3)$-face of the fixed simplex $A_0\dots A_n$ one of the $(n\minus 2)$-dimensional symmetry spaces according to \Thmref{thm:twice-flat n} passes through a common (finite or infinite) point $\Sq$.
\end{thm}

\begin{Proof}
A cross-polytope $\Ccal n$ is flexible if the simultaneous continuous rotations of the wings preserve all distances $\ol{B_iB_j}$ with $i\ne j$.
For twice-flat $\Ccal n$ this means by \Lemref{lem:4d-transmission} that for each $(i,j)$ the tangents $t_i, t_j$ of the halved angles of rotation $\phi_i, \phi_j$ satisfy \eqref{eq:tan_proportional_n} for any point $\Sq_{ij}$ in one of the $(n\minus 2)$-planes of symmetry of the pairs $\left( [\Kante_{ij}, A_i],\,[\Kante_{ij}, \Bq_i]\right)$ or $\left([\Kante_{ij}, A_j],\,[\Kante_{ij}, \Bq_j]\right)$.
\\[0.5mm]
(i) If there exists a common center $\Sq$ for all $(i,j)$, then for given $\phi_1$ the choice
\[   t_i:= \frac{d(\Sq,\Achse_1)}{d(\Sq,\Achse_i)}\,t_1 \quad \mbox{for} \ i=2,\dots,n
\]
guarantees flexibility since for every two $t_i,t_j$ eq.\ \eqref{eq:tan_proportional_n} holds true.

\smallskip\noindent
(ii) Conversely, if a given twice-flat cross-polytope $\Ccal n$ is flexible, then we use the vertices $A_1,\dots, A_n$ of the fixed simplex $\Rastn$ as base points of a projective coordinate frame in the projective extension of $\Rasthe$ and the center of the insphere of $\Rastn$ as its unit point.\footnote{
This choice is equivalent to the usage of signed Hessian normal forms for the $n$ bounding hyperplanes of the simplex $\Rastn$ in order to define the coordinates.}
Then, the homogeneous coordinates of any finite point in $\Rasthe$ are proportional to the signed distances from the facets of $\Rastn$.
\\ 
The transmissions between the rotations of any two wings $\Wing_i, \Wing_j$ of the given cross-polytope $\Ccal n$ define by \eqref{eq:tan_proportional_n} homogeneous coordinates $(t_1:\dots:t_n)$ of a unique point $\Sq\in\Rasthe$, and by \eqref{eq:transmission1} and \Lemref{lem:projection} the point $\Sq$ is located on one of the axes of local symmetry at each $(n\minus 3)$-face of $\Rastn$. 
\end{Proof}

Like in $\Raum E^3$, the motion of the moving simplex $\Gangn$ relative to $\Rastn$ defines a projective mapping between the circular trajectories of every two moving points $B_i, B_j$. 

\begin{cor}\label{cor:concurrent En} 
If in $\Rasthe$ at a flat pose $A_1\dots \Bq_n$ of a twice-flat cross-polytope $\Ccal n$ through each $(n\minus 3)$-face of the fixed simplex $\Rastn = A_1\dots A_n$ one of the $(n\minus 2)$-dimensional symmetry spaces in the sense of \Thmref{thm:twice-flat n} passes through a common point $\Sq$, then the same is true for all $(n\minus 3)$-faces in the flat pose $A_1\dots \Bq_n$ of $\Ccal n$.
\end{cor}

\begin{Proof}
By virtue of \Thmref{thm:fund_thm_En} the existence of $\Sq$ characterizes flexible twice-flat cross-polytopes.
However, for the flexibility it is not relevant which facet of $\Ccal n$ is supposed to be fixed.
The center $\Sq$ remains the same since we can proceed from $\Rastn$ to any other facet step by step by iteratively replacing one vertex $A_i$ by $B_i$.
\end{Proof}

Due to \Thmref{thm:fund_thm_En}, a flat pose of a flexible twice-flat cross-polytope $\Ccal n$ can be constructed as given below.
The three-dimensional version dates back to Bricard in \cite[p.~144]{Bricard1}.

\medskip\noindent
\textbf{Construction.}
For obtaining a flat pose of a flexible twice-flat cross-polytope in $\Raum E^n$, $n\ge 3$, 
\\[0.5mm]
1.\ choose within a hyperplane $\Rasthe$ of $\Raum E^n$ a simplex $\Rastn = A_1\dots A_n$ and a point $\Sq$ off from the bounding hyperplanes and from hyperplanes that bisect any dihedral angle of $\Rastn$.
\\[0.5mm]
2.\ For all $(i,j)$ with $i,j\in\{1,\dots,n\}$ and $i\ne j$, let $\Kante_{ij}$ denote the span of the $(n\minus 3)$-face of $\Rastn$ opposite to $A_iA_j$.
Then, reflect $\Achse_j = [\Kante_{ij},A_i]$ in the connection $[\Kante_{ij}, \Sq]$ to obtain the $(n\minus 2)$-dimensional space $g_{ij}$ passing through $\Bq_i$  (\Figref{fig:proj}).
\\[0.5mm]
3.\ Within the hyperplane $\Rasthe$, the point $\Bq_i$ is the intersection of the $(n\minus 2)$-planes $g_{ij}$ for $j\in\{1,\dots,n\} \setminus\{i\}$ (see the choices $(i,j)=(1,2)$ or $(2,1)$ in \Figref{fig:proj}).

\medskip\noindent
\begin{rem} The choice of $\Sq$ in a hyperplane that bisects the dihedral angle of $\Rastn$ at $\Kante_{ij}$ yields flat wings $\Wing_i$ and $\Wing_j$ (see \Figref{fig:proj}), which has been excluded from the beginning.
\end{rem}

\begin{thm}\label{thm:flat_flexible n} 
Let $A_1 \dots \Bq_n$ in the hyperplane $\Rasthe$ of $\Raum E^n$, $n\ge 4$, be a flat pose of a flexible twice-flat cross-polytopes $\Ccal n$ with a finite center $\Sq$ of local symmetries.
Then holds: 
\\[0.5mm]
{\em (i)} For each $i\in\{1,\dots,n\}$, the spans of all $(n\minus 2)$-faces with vertices different from $A_i$ and $\Bq_i$ contact the same hypersphere $\Scal_i$ centered at $\Sq$.
Consequently, each $(n\minus 2)$-face of the given flat pose of $\Ccal n$ has a span which is tangent to one of in total $n$ concentric hyperspheres $\Scal_1, \dots, \Scal_n$ in $\Rasthe$.
\\[0.5mm]
{\em (ii)} The $n$ dual quadrics in $\Rasthe$ of rank $2$ consisting of two bundels of hyperplanes with respective carriers $A_i$ and $\Bq_i$ for  $i = 1,\dots,n$ span a linear set that contains the rank-1 quadric of hyperplanes through $\Sq$ and the family of isotropic hyperplanes.
\\[0.5mm]
{\em (iii)} The midpoints of the diagonals $A_i\Bq_i$ for $i=1,\dots,n$ are coplanar with the center $\Sq$ of local symmetries.
\\[0.5mm]
{\em (iv)} If $n$ is even, then the second flat pose $\Bqq_1\dots \Bqq_n$ of the moving simplex $\Gangn$ is directly congruent to the first $\Bq_1\dots \Bq_n$ within $\Rasthe$.
Otherwise the congruence transformation in $\Rasthe$ from $\Gangnq$ to $\Gangnqq$ is orientation reversing.
\\[0.5mm]
{\em (v)} Let $\Sqq$ be the center of local symmetry in the second flat pose and $\Kante_{ij} = \Achse_i\cap\Achse_j$ be the span of any $(n\minus 3)$-face of the fixed simplex $\Rastn$.
Then, the connections $[\Kante_{ij},\Sq]$ and $[\Kante_{ij},\Sqq]$ as well as the facets $\Achse_i$ and $\Achse_j$ of $\Rastn$ share the angle bisecting hyperplanes within $\Rasthe$.
In this sense, the points $\Sq$ and $\Sqq$ are isogonal w.r.t.\ $A_1\dots A_n$.
\end{thm}

\begin{Proof}
(i) By virtue of \Corref{cor:concurrent En}, the local symmetry of the given flat pose of $\Ccal n$ implies within the $(n\minus 1)$-dimensional $\Rasthe$ that the connection of $\Sq$ with any $(n\minus 3)$-face of $\Ccal n$ is a hyperplane of local symmetry.
This means for example, that the reflection in $[\Sq, X_3\dots X_n]$ with $X_k\in\{A_k,B_k\}$ exchanges $[A_2X_3\dots X_n]$ with $[B_2X_3\dots X_n]$.
Consequently, these two $(n\minus 2)$-spaces have equal distances to $\Sq$.
Iteration shows that $[B_2A_3\dots X_n]$ and $[B_2B_3\dots X_n]$ have equal distances to $\Sq$ due to their symmetry w.r.t.\ $[\Sq,B_2X_4\dots X_n]$, and so on.
Thus, we obtain $2^{n-1}$ different $(n\minus 2)$-faces of the flat cross-polytope in $\Rasthe$ with spans at equal distances to $\Sq$.
These $(n\minus 2)$-faces are charactered by containing no vertex with the subscript $1$. 
\\
Statement (i) is apparently the $n$-dimensional generalization of \Figref{fig:Fig_1}(a); the families of $2^{n-1}$ $(n\minus 2)$-faces of $\Ccal n$ tangent to the same hypersphere are the counterparts of the quadrangles on octahedra.
  
\smallskip\noindent
(ii) We prove that the hyperspheres $\Scal_1,\dots,\Scal_n$ belong to the addressed linear system of quadrics in $\Rasthe$.
As a consequence, also the dual pencil of concentric hyperspheres is included und consequently the mentioned rank-1 quadric and the family of isotropic hyperplanes in $\Rasthe$.

\noindent
It is sufficient to show that $\Scal_1$ is contained in the subsystem spanned by the rank-2 quadrics with carriers $(A_2,\Bq_2)$, \dots, $(A_n,\Bq_n)$, which are the vertices of a full-dimensional cross-polytope in $\Rasthe$:
\\
Each quadric of this subsystem in $\Rasthe$ contains the hyperplanes $[X_2\dots X_n]$ since they are common to all spanning rank-2 quadrics.
Thus, the subsystem is a subset of the dual quadrics containing the $2^{n-1}$ hyperplanes.
\\
Conversely, the dual quadrics in $\Rasthe$ through these hyperplanes form a linear system, since the entries in the related symmetric coefficient matrices are the solutions of a linear system of equations.
For each quadric $\Qcal$ out of this linear system holds for each  $j\in \{2,\dots,n\}$:
If a hyperplane of $\Qcal$ passes through one point of the pairs $(A_i,\Bq_i)$ for all $i\in\{2,\dots,n\}$ with $i\ne j$, then it also passes through $A_j$ or $\Bq_j$.
Passing through $A_i$ or $\Bq_i$ is for hyperplanes equivalent to be contained in the rank-2 quadric with carriers $(A_i,\Bq_i)$.
Thus, $\Qcal$ belongs to the $(n\minus 2)$-dimensional subsystem.
With other words, each full-dimensional cross-polytope $\Ccal{n-1}$ in $\Rasthe$ defines an $(n\minus 2)$-dimensional linear system of quadrics that contain the spans of all facets of $\Ccal{n-1}$.
\\
Due to (i), the hypersphere $\Scal_1$ with center $\Sq$ contacts all extended facets of the cross-polytope with vertices $A_2\dots A_n\Bq_2\dots \Bq_n$ and is therefore contained in the subsystem defined above.
Similar properties of the other hyperspheres $\Scal_2\dots\Scal_n$ confirm the statement (ii).

\smallskip\noindent
(iii) In $\Rasthe$, the dual quadrics contained in a $(n\minus 1)$-parametric linear system have their centers located in a hyperplane. 

\smallskip\noindent
(iv) Below there is a sequence of $n\plus 1$ facets of $\Ccal n$ beginning with $\Gangn$ and ending with $\Rastn$.
Every two consecutive facets are neighboring, i.e., they share an $(n\minus 2)$-face that spans the hinge between the facets.
\def\abst{\hspace{2.0mm}}
\def\abstb{\hspace{2.8mm}}
{\small
\begin{equation}\label{eq:sequence}
 \begin{array}{ll} 
   \mathrm{facets:\!}  &B_1B_2\dots B_n, \abst A_1B_2\dots B_n, \abst 
     A_1A_2B_3\dots B_n, \abstb \dots, \abstb A_1\dots A_{n-1}B_n, \abst 
     A_1A_2\dots A_n
   \\[1.0mm]
   \mathrm{hinges:}   &\hspace{10.5mm} [B_2\dots B_n], \abstb [A_1B_3\dots B_n], 
   \abstb [A_1A_2B_4\dots B_n], \abst \dots,\abst [A_1A_2\dots A_{n-1}].
 \end{array}
\end{equation}
}We use this sequence to explain how the first flat pose $\Gangnq = \Bq_1\dots \Bq_n$ of $\Gangn$ can be transferred in the second flat pose $\Gangnqq = \Bqq_1\dots \Bqq_n$ by $n$ consecutive halfturns:
\\
A first halfturn is carried out about the hinge $[\Bq_2\dots\Bq_n]$ between $\Gangnq$ and $A_1\Bq_2\Bq_3\dots \Bq_n$ in order to achieve the correct relative position of these two facets according to the second flat pose.
Afterwards, the moving simplex is attached to $A_1\Bq_2\dots \Bq_n$ and performs the halfturn about the hinge $[A_1\Bq_3\dots\Bq_n]$ between $A_1\Bq_2\dots \Bq_n$ and $A_1A_2\Bq_3\dots \Bq_n$ to achieve the correct positions of the first two facets relative to the third one according to the other flat pose.
Afterwards, $\Gangn$ is attached to $A_1A_2\Bq_3\dots \Bq_n$ during the halfturn about the hinge $[A_1A_2\Bq_4\dots\Bq_n]$ with $A_1A_2A_3\Bq_4\dots \Bq_n$, and so on.
The final halfturn uses the hinge $\Achse_n = [A_1\dots A_{n-1}]$ between $A_1\dots A_{n-1}\Bq_n = \Wing_n$ and $A_1\dots A_n = \Rastn$.
\\
Within the hyperplane $\Rasthe$, which contains both flat poses, the halfturns act like reflections in $(n\minus 2)$-dimensional spaces, i.e., hyperplanes of $\Rasthe$, and each of them reverses the orientation. 
In \Thmref{thm:spatial_n},\,(iii) we show that there exist $n!$ different sequences with the same property like that in \eqref{eq:sequence}.

\smallskip\noindent
(v) This is a consequence of the local symmetries of the two flat poses $A_1\dots\Bq_n$ and $A_1\dots\Bqq_n$ and the symmetry of $\Bq_i$ and $\Bqq_i$ w.r.t.\ the hinge $h_i$ (note \Figref{fig:proj}): 
The product of the reflections in $[v_{12},\Sq]$ and $h_1$ takes $[v_{12},A_1]$ via $[v_{12},\Bq_1]$ to $[v_{12},\Bqq_1]$.
\end{Proof}

The following theorem is an extension of the statement (ii) in \Thmref{thm:flat_flexible n}.
Here we use the symbol $\Vkt M_{XY}$ for the symmetric coefficient matrix of the rank-2 dual quadric of hyperplanes passing through the points $X$ or $Y$ in the hyperplane $\Rasthe$ of $\Raum E^n$.
This means, if $(1:x_1:\dots:x_{n-1})$ and $(1:y_1:\dots:y_{n-1})$ are homogeneous Cartesian coordinates of $X$ and $Y$, that 
\begin{equation}\label{eq:matrix}
  2\,\Vkt M_{XY} := (1\, x_1 \dots x_{n-1})^\top(1\, y_1\dots y_{n-1}) +  
     (1\, y_1 \dots y_{n-1})^\top(1\, x_1\dots x_{n-1}).
\end{equation}

\begin{thm}\label{thm:characterization} 
The cross-polytope $\Ccal n$ given by the flat pose $A_1\dots \Bq_n$ is twice-flat and flexible in $\Raum E^n$, $n\ge 3$, if and only if the set of linear combinations of the $(n\mal n)$-matrices 
\[ \Vkt M_{A_1\Bq_1}, \ \Vkt M_{A_2\Bq_2}, \ \dots, \ \Vkt M_{A_n\Bq_n}
\]
contains a rank-1 matrix of type $\Vkt M_{\Sq\Sq}$ and the diagonal matrix $\mathrm{diag}(0,1,\dots,1)$.
\end{thm}

\begin{Proof}
The dual quadrics of the linear system addressed in \Thmref{thm:flat_flexible n}\,(ii) have coefficient matrices that are linear combinations of $\Vkt M_{A_1\Bq_1}, \,\dots,\, \Vkt M_{A_n\Bq_n}$.
In order to prove the converse of \Thmref{thm:flat_flexible n}\,(ii), we choose any $(n\minus 3)$-face of $\Rastn$ and confirm the local symmetry within $\Rasthe$.
It means no restriction of generality to assume that this is the face $A_3\dots A_n$ with the span $v_{12}$.
\\
From the viewpoint of Projective Geometry, for each dual quadric out from the linear system, the hyperplanes passing through the subspace $v_{12}$ form again a quadric in the bundle except that the whole bundle is part of the quadric.
The latter holds for all rank-2 quadrics with carriers $(A_i,\Bq_i)$ for $i\ge 3$.
Thus, the restriction of our linear system to the bundle with carrier $v_{12}$ is spanned by the restrictions of the two rank-2 quadrics with carriers $(A_1,\Bq_1)$ and $(A_2,\Bq_2)$.
The orthogonal view in \Figref{fig:proj} with the point $v_{12}^s$ as image of $v_{12}$ visualizes the restriction of the linear system as the dual pencil of conics tangent to the sides of the quadrangle $A_1^s A_2^s \Bqn_1 \Bqn_2$.
Due to Desargues's involution theorem (see, e.g., \cite[Sect.~7.4]{Conics}), the tangents from $v_{12}^s$ to the conics form an involution where the connection with $\Sqn$ is one fixed line.
Since in the orthogonal view also the isotropic lines through $v_{12}^s$ are a pair of corresponding lines, the two fixed lines are orthogonal and the involution with the particular pairs $([v_{12}^s,A_1^s],\,[v_{12}^s,\Bqn_1])$ and $([v_{12}^s,A_2^s],\,[v_{12}^s,\Bqn_2])$ is symmetric w.r.t.\ $[v_{12}^s,\Sqn]$.
\\
The same holds for all $(n\minus 3)$-faces of $\Rastn$, and this confirms the converse direction of the stated equivalence.  
\end{Proof}

Note that the first row of the matrices in \Thmref{thm:characterization} confirms again the statement (iii) in \Thmref{thm:flat_flexible n}.

\subsection{Spatial poses of flexible twice-flat cross-polytopes}

\begin{thm}\label{thm:spatial_n} 
Let $A_1\dots B_n$ be a spatial pose of a flexible twice-flat cross-polytope $\Ccal n$ in $\Raum E^n$, $n\ge 4$, with a finite center $\Sq$ of local symmetries in a flat pose.
Then holds for each permutation $(i_1,\dots,i_n)$ of $(1,\dots,n)$:
\\[0.5mm]
{\em (i)} The dihedral angle along the $(n\minus 2)$-dimensional hinge $[X_{i_2}\dots X_{i_n}]$ between the neighboring facets $A_{i_1} X_{i_2} \dots X_{i_n}$ and $B_{i_1} X_{i_2} \dots X_{i_n}$ of $\Ccal n$, where $X_k\in\{A_k,B_k\}$, is either congruent or supplementary to the angle $\phi_{i_1}$ between the wing $\Wing_{i_1}$ and the fixed simplex $\Rastn$.
Congruence holds if and only if at the initial flat pose the corresponding dihedral angle equals $0$.
Thus, the $2^{n-1}\!\cdot n$ dihedral angles of $\Ccal n$ can be subdivided into $n$ classes of congruent or supplementary angles. 
\\[0.5mm]
{\em (ii)} The one-parameter motion of the simplex $\Gangn$ relative to $\Rastn$ can be expressed as a product of $n$ rotations through respective angles $\phi_1, \dots, \phi_n$ about fixed hinges spanned by $(n\minus 2)$-faces of the flat pose $A_1\dots \Bq_n$ in the hyperplane $\Rasthe$. 
There exist $n\klzwi!$ different decompositions of this kind.
The induced motion of a facet with $m$ different $B$-vertices $B_{i_1},\dots,B_{i_m}$, $m<n$, is the product of $m$ rotations through respective angles 
$\phi_{i_1},\dots,\phi_{i_m}$ about hinges in $\Rasthe$.
\\[0.5mm]
{\em (iii)} Among the four $(n\minus 2)$-faces  of $\Ccal n$ passing through the $(n\minus 3)$-face $X_{i_3}\dots X_{i_n}$, the connections with $A_{i_1}$ and $B_{i_1}$ as well as that with $A_{i_2}$ and $B_{i_2}$ have a common angle bisecting hyperplane $\mu$. 
This hyperplane $\mu$ through $\Kante:= [X_{i_3}\dots X_{i_n}]$ bisects also the angle between the pairs of opposite hyperplanes $\left([A_{i_1}A_{i_2},\Kante],\  [B_{i_1}B_{i_2},\Kante]\right)$ and of $\left([A_{i_1}B_{i_2},\Kante],\ [B_{i_1}A_{i_2},\Kante]\right)$, which extend facets of $\Ccal n$. 
\\[0.5mm]
{\em (iv)} All these angle-bisecting hyperplanes $\mu$ of $\Ccal n$ share a line $s$ in $\Raum E^n$.
All hinges of type $[X_{i_2}\dots X_{i_n}]$ enclose the same angle $\psi_{i_1}$ with $s$.
\\[0.5mm]
{\em (v)} All facets of $\Ccal n$ with an even number of $A$-vertices span hyperplanes that share a point $C\in s$.
Similarly, all hyperplanes spanned by facets with an odd number of $A$-points pass through another common point $D\in s$.
\\[0.5mm]
{\em (vi)} All $2^n$ facets of $\Ccal n$ are placed in hyperplanes that contact a hypersphere centered on $s$.
The midpoints of the diagonals $A_iB_i$ belong to a hyperplane passing through $s$.  
\end{thm}

\begin{Proof}
(i) The four facets of $\Ccal n$ through the $(n\minus 3)$-face $X_{i_3} \dots X_{i_n}$ form an isogonal pyramid where opposite dihedral angles are congruent of supplementary (\Figref{fig:proj_b}).
Hence, the dihedral angle between the neighboring faces $A_{i_1}A_{i_2}X_{i_3} \dots X_{i_n}$ and $B_{i_1}A_{i_2}X_{i_3} \dots X_{i_n}$ along the hinge
$[A_{i_2}X_{i_3} \dots X_{i_n}]$ is congruent or supplementary to that between $B_{i_1}B_{i_2}X_{i_3} \dots X_{i_n}$ and $A_{i_1}B_{i_2}X_{i_3} \dots X_{i_n}$ along the hinge $[B_{i_2}X_{i_3} \dots X_{i_n}]$; in both involved hinges the vertices with index $i_1$ are missing. 
In other words, a dihedral angle at $\Ccal n$ remains congruent or becomes supplementary if at a hinge one vertex is replaced by its opposite.
Iteration confirms the stated property.

\smallskip\noindent
(ii) We refer to the sequence of neighboring facets in \eqref{eq:sequence}, $\Gangn = B_1\dots B_n$, $A_1B_2\dots B_n$, $A_1A_2B_3\dots B_n$, $\dots$, $ A_1\dots A_n = \Rastn$ with $n$ hinges between as $[B_2\dots B_n]$, $[A_1B_3\dots B_n]$, $\dots$, $[A_1\dots A_{n-1}] = h_n$.
According to (i), at the spatial pose of $\Ccal n$ the dihedral angles along these hinges are congruent or supplementary to $\phi_1, \phi_2, \dots, \phi_n$, respectively. 
Hence, we obtain the general pose stepwise by appropriate rotations about the hinges in the flat pose in order to establish these angles between the neighboring faces.
\\
Here we face the problem that the angles $\phi_1,\dots,\phi_n$ of rotations about  the hinges in $\Rasthe$ are oriented.
How to orientate an $(n\minus 2)$-dimensional hinge?
\\
We follow the example of octahedra (see \Figref{fig:denotation}) and use the orthogonal projection along an $(n\minus 3)$-dimensional subspace of the hinge $h$ similar to that in \Figref{fig:proj}.
Under the condition that the projection yields an orientation-preserving view of a predefined side of $\Rasthe$, we orientate the view $h^s$ of the hinge $h$ such that the interior of $\Rastn$ lies left. 
Now we define that, when looking in this direction, a rotation through a positive angle is carried out counter-clockwise.
\\
The sequence of facets in \eqref{eq:sequence} leads to the natural order $\phi_1,\phi_2,\dots,\phi_n$ of the angles.
And this order corresponds to that of the consecutive rotations.
However, we can also start by predefining any order of the angles.
Then, conversely, this determines the sequence of the facets according to (i) and consequently the sequence of the hinges and the rotations.
This yields $n\klzwi !$ different representations.
In a similar way we proceed when we represent the movement of any other facet with $B$-vertices $B_{i_1},\dots,B_{i_m}$, $m<n$, against $\Rastn$.
\\
The 3d-versions of these representations are equivalent to the different decompositions of the motion polynomial of a type-3 Bricard octahedron as discussed in \cite{Hege}.

\smallskip\noindent
(iii) This follows from the symmetry of the isogonal double-pyramid with the $(n\minus 3)$-dimensional apex space $\Kante = [X_{i_3}\dots X_{i_n}]$ and `edges' through the vertices $A_{i_1}, A_{i_2},B_{i_1},B_{i_2}$ (note the orthogonal view in \Figref{fig:proj_b} with $(i_1,i_2) = (i,j)$).
Obviously, this property of $\Ccal n$ is the spatial analogue of the local symmetry of flat poses in $\Rasthe$ as stated in \Thmref{thm:twice-flat n}:
From the viewpoint of $\Raum E^n$, all angle bisecting hyperplanes of the connections of $\Kante$ with $(A_{i_1},\Bq_{i_1})$ and with $(A_{i_2},\Bq_{i_2})$ share the perpendicular $s'$ to $\Rasthe$ through the center $\Sq$ (note \Figref{fig:proj_b}).

\smallskip\noindent
(iv) The symmetries of the double-pyramids assign to each $(n\minus 3)$-face of $\Ccal n$ a unique hyperplane $\mu$ of symmetry between the spans of the pyramid's opposite facets.
In addition, this angle-bisecting plane $\mu$ is already uniquely defined by one of the two pairs of opposite hyperplanes, as---by virtue of the flexion---it can be continuously transferred into the corresponding hyperplane through $s'$ in the flat pose.
Opposite facets of such a pyramid are those which differ at exactly two vertices while the number of $A$-points either changes by $2$ or remains the same.
Of course, the same holds for the number of $B$-points.
Why do all $2^{n-3} n(n\minus 1)$ hyperplanes $\mu$, that are assigned to the $(n\minus 3)$-faces of $\Ccal n$, share a line $s\klzwi$?
\\[0.5mm]
At the begin, we focus on the $(n\minus 3)$-faces of $\Rastn = A_1\dots A_n$ and refer to the previously introduced notation $\Kante_{ij}$ for the span of all $A$-vertices different from $A_i$ and $A_j$.
The assigned bisecting hyperplane will be denoted by $\mu_{ij}$, and we define the line $s$ as the intersection of $n\minus 1$ hyperplanes, namely
\begin{equation}\label{eq:s}
  s: = \mu_{12}\cap\mu_{23}\cap\mu_{34}\cap\mu_{45}\cap\dots\cap\mu_{n-1\,n}.
\end{equation}
Hence, each point $S\in s$ has equal distances to the pairs of hyperplanes $\left([A_iA_{i+1},\Kante_{i\Zw i+1}],\right.$ $\left.[B_iB_{i+1}, \Kante_{i\Zw i+1}]\right)$ and $\left([A_iB_{i+1}, \Kante_{i\Zw i+1}], [B_iA_{i+1}, \Kante_{i\Zw i+1},]\right)$ for $i=1,\dots,n\minus 1$.
This defines an equivalence relation between facets for which we introduce the symbol `$\sim$'.
For consecutive $\Kante_{i\,i+1}$ follow by \eqref{eq:s} the equivalences
\begin{equation}\label{eq:s2}
 \begin{array}{l}
    \Rastn \sim B_1B_2A_3\ldots \sim A_1B_2B_3A_4\ldots 
     \sim A_1A_2B_3B_4A_5\ldots \sim \dots \sim A_1\ldots A_{n-2}B_{n-1}B_n,
    \\[1mm]
    A_1B_2A_3\dots A_n \sim B_1A_2A_3 \dots A_n , \
     A_1A_2B_3A_4\dots A_n \sim A_1B_2A_3\dots A_n, 
    \\[0.5mm]
    \dots, \ A_1\dots A_{n-1}B_n \sim A_1\dots A_{n-2}B_{n-1}A_n.
  \end{array} 
\end{equation}   
Our goal is to extend the list of the involved $(n\minus 3)$-faces step by step so that finally all $2^{n-3}\cdot n(n\minus 1)$ apex-spaces of pyramids on $\Ccal n$ are involved. 

\noindent 
1.\ The first line in \eqref{eq:s2} shows that opposite facets through the $(n\minus 3)$-dimensional apex-spaces
\begin{equation}\label{eq:s3}
 [B_2A_4A_5\dots], \  [A_1B_3A_5\dots], \ \dots, \ [A_1\dots A_{n-3}B_{n-1}]
\end{equation}   
are equivalent.
The first equivalence in the second line of \eqref{eq:s2} along with the first $(n\minus 3)$-face in \eqref{eq:s3} yield
\begin{equation}\label{eq:s4}
   B_1A_2A_3A_4\ldots \sim A_1B_2A_3A_4\ldots \sim B_1B_2B_3A_4\ldots 
\end{equation}   
as opposite facets w.r.t.\ $[B_1A_4\dots]$.
Together with the first equivalence in \eqref{eq:s2} follows
\begin{equation}\label{eq:s5}
   \Rastn \sim B_1B_2A_3A_4\ldots \sim B_1A_2B_3A_4\ldots
\end{equation}   
as opposite facets w.r.t.\ $[A_2A_4A_5\dots] = \Kante_{13}$.
\\
When (under $n>4$) we iteratively increase the involved subscripts by $1$ as long as none of them exceeds $n$, then we obtain a sequence of $(n\minus 3)$-dimensional apex-spaces $\Kante_{13}, \, \Kante_{24}, \, \dots, \, \Kante_{n-2\Zw n}$ with associated hyperplanes $\mu_{13}, \mu_{24}, \mu_{35}, \dots, \mu_{n-2\Zw n}$ 	passing through $s$.

\noindent 
2.\ The definition of $s$ remains valid if we exchange in \eqref{eq:s2} the indices $3$ and $4$.
This implies that the planes $\mu_{14}, \mu_{25}, \mu_{36}, \dots, \mu_{n-3\Zw n}$ pass through $s$, too.
We repeat by exchanging $4$ and $5$ and get $s\subset \mu_{15}, \mu_{26}, \dots, \mu_{n-4\Zw n}$.
Iteration yields finally $s\subset\mu_{1n}$ so that all planes of symmetry $\mu_{ij}$ assigned to the $(n\minus 3)$-faces of $\Rastn$ are proved to contain the line $s$ as defined in \eqref{eq:s}.

\noindent 
3.\ We continue with other facets of $\Ccal n$.
Each facet $\Wing_{ij} = A_1\dots B_i\dots B_j\dots A_n$ that contains two $B$-points is opposite to $\Rastn$ w.r.t.\ $\Kante_{ij}$.
Each facet $\Wing_i = A_1\dots B_i\dots A_n$ with exactly one $B$-point is opposite to $\Wing_1$ w.r.t.\ $\Kante_{1i}$.
Thus, all facets with one or two $B$-points can be inserted in our two equivalence classes.
Moreover, the $(n\minus 3)$-faces of $\Wing_i$ are either part of $\Rastn$ or they contain one $B$-point.
If it contains $B_i$ but not $A_j$ and $A_k$, $j\ne k$, then at the corresponding pyramid the equivalent facets $\dots B_i\dots A_j\dots B_k\dots$ and  $\dots B_i\dots B_j\dots A_k\dots$ with two $B$-points are opposite so that the plane of symmetry assigned to this $(n\minus 3)$-face of $\Wing_i$ passes again through $s$.
\\
Hence, we can replace $\Rastn$ by $\Wing_{i}$ and start the procedure again.
As a result, all facets containing $B_i$ and two further $B$-points are equivalent to $\Wing_i$.
Moreover, we can prove like before that all $(n\minus 3)$-faces of $\Wing_{ij}$ which contain two $B$-points are the meet of two opposite and equivalent facets as they contain three $B$-points.
Thus, for all apex-spaces of $\Wing_{ij}$ the assigned hyperplanes $\mu$ pass through $s$.
Now, we iterate and replace $\Rastn$ by $\Wing_{ij}$.
This allows to extend our equivalence classes to facets of $\Ccal n$ with four $B$-points, and so on, until we reach $\Gangn$ with $n$ $B$-points. 

\noindent
4.\ Let $\mu$ be the plane of symmetry assigned to the $(n\minus 3)$-face $X_{i_3}\dots X_{i_n}$.
The reflection in $\mu$ exchanges opposite hinges of the double-pyramid and preserves $s\subset\mu$ (see \Figref{fig:proj_b}).
Thus, the angles between $s$ and opposite hinges $[A_{i_2}X_{i_3}\dots X_{i_n}]$ and $[B_{i_2}X_{i_3}\dots X_{i_n}]$ are congruent.
By iterated transition to opposite hinges one obtains only hinges that contain neither $A_{i_1}$ or $B_{i_1}$.
These are exactly the hinges of $\Ccal n$ with dihedral angles congruent or supplementary to $\phi_{i_1}$ as mentioned in (i).

\smallskip\noindent
(v) Each point $S\in s$ has equal distances to the spans of equivalent facets.
According to the two equivalence classes, the distances define two linear functions on $s$. 
The zeros of these functions show that the hyperplanes of each class have a point on $s$ in common, and these two points $C,D\in s$ must be different since otherwise all hinges of $\Ccal n$ would share a point.
We cannot exclude that one of these points lies at infinity.

\smallskip\noindent
(vi) The fixed simplex $\Rastn$ and the wing $\Wing_1$ are representatives of the two equivalence classes.
An orthogonal projection of $\Raum E^n$ along the hinge $\Achse_1$ into a plane shows $\Rasthe$ and $[\Wing_1]$ as two lines enclosing the angle $\phi_1$.
In the image plane, the two angle bisectors intersect the image line of $s$ at two points $S_1,S_2$ in harmonic position w.r.t.\ $C$ and $D$.
Since distances of points $S\in s$ to $\Rasthe$ and $[\Wing_1]$ appear in true size, each $S_i$ has equal distances to both and, consequently, to all other hyperplanes, provided that $S_i$ is finite.
\\
The proof of \Thmref{thm:flat_flexible n},\,(ii) revealed that cross-polytopes in $\Rastn$ define an $(n\minus 2)$-parametric linear system of dual quadrics containing all hyperplanes spanned by the cross-polytope's facets.
The same holds for $\Ccal n$ in $\Raum E^n$.
The $(n\minus 1)$-parametric linear system is spanned by the rank-2 quadrics with carriers $(A_i,B_i)$ for $i=1,\dots,n$ and contains also the two hyperspheres with centers $S_1,S_2\in s$.
All dual quadrics of this linear system have their centers in a common hyperplane which must pass through $s$.
Continuity guarantees this also when one $S_i$ tends to infinity.
\end{Proof}

\begin{rem}\label{rem:Grundriss_n}
The comment in Remark~\ref{rem:Grundriss} can be generalized to $\Raum E^n$ due to \Thmref{thm:spatial_n}\,(iii):
An orthogonal projection of any spatial pose of a flexible twice-flat cross-polytope $\Ccal n$ in direction of the line $s$ yields an $(n\minus 1)$-dimensional image that has the properties of flat poses as listed in \Thmref{thm:flat_flexible n}.
\end{rem}


\end{document}